\documentclass[11pt]{article}
\usepackage[english]{babel}
\usepackage[latin1]{inputenc}
\usepackage{amsmath}
\usepackage{amssymb}
\usepackage{amsfonts}
\usepackage{amsthm}

\newcommand{\p}{\mathfrak{p}}

\newcommand{\C}{\mathbb{C}}
\newcommand{\N}{\mathbb{N}}

\newcommand{\R}{\mathbb{R}}

\newcommand{\boC}{\mathcal{C}}
\newcommand{\boE}{\mathcal{E}}

\newcommand{\boL}{\mathcal{L}}

\newcommand{\boP}{\mathcal{P}}
\newcommand{\boR}{\mathcal{R}}
\newcommand{\boS}{\mathcal{S}}

\newcommand{\eps}{\varepsilon}
\newcommand{\ch}{{\rm ch}}

\renewcommand{\div}{{\rm div}}

\renewcommand{\Im}{{\rm Im}}
\newcommand{\on}{\ {\rm on} \ }
\renewcommand{\Re}{{\rm Re}}

\newtheorem{cor}{Corollary}
\newtheorem{lemma}{Lemma}
\newtheorem{prop}{Proposition}
\newtheorem{step}{Step}
\newtheorem{theorem}{Theorem}
\theoremstyle{definition}
\newtheorem*{acknowledgement}{Acknowledgements}
\newtheorem{remark}{Remark}

\begin{document}

\title{On the KP I transonic limit of two-dimensional Gross-Pitaevskii travelling waves}
\author{
\renewcommand{\thefootnote}{\arabic{footnote}}
Fabrice B\'ethuel \footnotemark[1], Philippe Gravejat \footnotemark[2], Jean-Claude Saut \footnotemark[3]}
\footnotetext[1]{Laboratoire Jacques-Louis Lions, Universit\'e Pierre et Marie Curie, Bo\^ite Courrier 187, 75252 Paris Cedex 05, France. E-mail: bethuel@ann.jussieu.fr}
\footnotetext[2]{Centre de Recherche en Math\'ematiques de la D\'ecision, Universit\'e Paris Dauphine, Place du Mar\'echal De Lattre De Tassigny, 75775 Paris Cedex 16, France. E-mail: gravejat@ceremade.dauphine.fr}
\footnotetext[3]{Laboratoire de Math\'ematiques, Universit\'e Paris Sud, B\^atiment 425, 91405 Orsay Cedex, France. E-mail: Jean-Claude.Saut@math.u-psud.fr}
\maketitle

\begin{abstract}
We provide a rigorous mathematical derivation of the convergence in the long-wave transonic limit of the minimizing travelling waves for the two-dimensional Gross-Pitaevskii equation towards ground states for the Kadomtsev-Petviashvili equation \eqref{KP}.
\end{abstract}

\section{Introduction}

\subsection{Statement of the results}

The Gross-Pitaevskii equation
\renewcommand{\theequation}{GP}
\begin{equation}
\label{GP}
i \partial_t \Psi = \Delta \Psi + \Psi (1 - |\Psi|^2) \on \R^N \times \R,
\end{equation}
appears as a relevant model in various areas of physics: Bose-Einstein condensation, fluid mechanics (see e.g. \cite{GinzPit1,Pitaevs1,Gross1,Coste1}), nonlinear optics (see e.g. \cite{KivsLut1})... At least on a formal level, this equation is hamiltonian, with a conserved Hamiltonian given by the Ginzburg-Landau energy,
\renewcommand{\theequation}{\arabic{equation}}
\setcounter{equation}{0}
\begin{equation}
\label{GLE}
E(\Psi) = \frac{1}{2} \int_{\R^N} |\nabla \Psi|^2 + \frac{1}{4} \int_{\R^N} (1 - |\Psi|^2)^2 \equiv \int_{\R^N} e(\Psi).
\end{equation}
Note that the boundedness of the Ginzburg-Landau energy implies that in some sense,
$$\vert \Psi(x, \cdot) \vert \to 1, \ {\rm as} \ \vert x \vert \to + \infty.$$
As a matter of fact, this condition provides a richer dynamics than in the case of null condition at infinity which is essentially governed by dispersion and scattering. In particular, \eqref{GP} has nontrivial coherent localized structures called travelling waves.

The existence of finite energy travelling waves was addressed and established in several papers (see \cite{IordSmi1,JoneRob1,JonPuRo1,BethSau1,BetOrSm1,Chiron1,BetGrSa1}). Travelling waves are special solutions to \eqref{GP} of the form
$$\Psi(x,t) = u(x_1 - c t, x_\perp), \ x_\perp = (x_2, \ldots, x_N).$$
They are supposed to play an important role in the full dynamics of \eqref{GP}. The equation for the profile $u$ is given by
\renewcommand{\theequation}{TWc}
\begin{equation}
\label{TWc}
i c \partial_1 u + \Delta u + u (1 - |u|^2) = 0.
\end{equation}
The parameter $c \in \R$ corresponds to the speed of the travelling waves. We may restrict to the case $c \geq 0$. Indeed, when $u$ is a travelling wave of speed $c$, the map $\overline{u}$ obtained by complex conjugation is a travelling wave of speed $- c$.

The existence of solutions to \eqref{TWc} was obtained in the above quoted papers through variational arguments, namely minimization under constraints \cite{BetOrSm1,BetGrSa1}, or mountain-pass theorems \cite{BethSau1,Chiron1}. In dimensions two and three, a full branch of solutions is constructed in \cite{BetGrSa1} minimizing the Ginzburg-Landau energy $E$ under fixed momentum $p$. In this context, the momentum is defined by
\renewcommand{\theequation}{\arabic{equation}}
\setcounter{equation}{1}
\begin{equation}
\label{VectP}
p(u) = \frac{1}{2} \int_{\R^N} \langle i \partial_1 u \ , u - 1 \rangle.
\end{equation}
This integral quantity is also formally conserved by \eqref{GP}. A notable difficulty in the variational approach is to give a meaning to the momentum in the space of maps of finite Ginzburg-Landau energy (see e.g. \cite{BetGrSa2,BeGrSaS1}). However, the momentum is well-defined for finite energy travelling wave solutions. Indeed, it is proved in \cite{Graveja3} that they belong to the space $W(\R^N)$, defined as
$$W(\R^N) = \{ 1 \} + V(\R^N),$$
where we have set
$$V(\R^N) = \big\{ v: \R^N \mapsto \C, \ {\rm s.t.} \ (\nabla v, \Re(v)) \in L^2(\R^N)^2, \Im(v) \in L^4(\R^N), \ {\rm and} \ \nabla \Re(v) \in L^\frac{4}{3}(\R^N) \big\}.$$
Separating real and imaginary parts, a direct computation shows that the quantity $\langle i \partial_1 v, v - 1 \rangle$ is integrable for any function $v \in W(\R^N)$, so that the momentum of travelling wave solutions is well-defined.

The main focus of this paper is a qualitative description of small Ginzburg-Landau energy solutions in the two-dimensional case. Such solutions are known to exist in view of the following result.

\begin{theorem}[\cite{BetGrSa1}]
\label{dim2}
i) Let $\p > 0$. There exists a non-constant finite energy solution $u_\p \in W(\R^2)$ to \eqref{TWc}, with
$0 < c = c(u_\p) < \sqrt{2}$, and
$$p(u_\p) \equiv \frac{1}{2} \int_{\R^2} \langle i \partial_1 u_\p, u_\p - 1 \rangle = \p,$$
such that $u_\p$ is solution to the minimization problem
$$E(u_\p) = E_{\min}(\p) = \inf \big\{ E(v), v \in W(\R^2), p(v) = \p \big\}.$$
ii) There exist some positive constants $K_0$, $K_1$ and $\boS_{KP}$, not depending on $\p$, such that we have the asymptotic behaviours
\begin{equation}
\label{estimE2}
0 < \frac{48 \sqrt{2}}{\boS_{KP}^2} \p^3 - K_0 \p^4 \leq \sqrt{2} \p - E(u_\p) \leq K_1 \p^3,
\end{equation}
for any $\p$ sufficiently small.
\end{theorem}

A more precise definition of the constant $\boS_{KP}$ will be provided in the course of our discussion of the Kadomtsev-Petviashvili equation \eqref{KP}. It should be noticed that we have, in view of \eqref{estimE2},
$$E(u_\p) \sim \sqrt{2} \p,$$
for small values of the momentum $\p$, so that Theorem \ref{dim2} provides a branch of travelling wave solutions with arbitrary small energy. Our aim is to describe the asymptotic behaviour, as $\p \to 0$, of the solutions $u_\p$ constructed above.

We recall that, in view of \cite{BethSau1,Graveja2,Graveja4}, any finite energy travelling waves are subsonic in dimension two, i.e. any non-constant finite energy solution $v$ to \eqref{TWc} satisfies
\begin{equation}
\label{subsonicwave}
0 < |c(v)| < \sqrt{2}.
\end{equation}
The speed $\sqrt{2}$ corresponds to the speed of sound waves at infinity around the constant solution $\Psi = 1$ to \eqref{GP}. Moreover, the quantity
$$\varepsilon(v) = \sqrt{2 - c(v)^2}$$
is related to the energy $E(v)$ and the uniform norm of $1 - |v|$ as follows.

\begin{prop}[\cite{BetGrSa1}]
\label{firstprop}
Let $v$ be a non-constant finite energy solution to \eqref{TWc} on $\R^2$. Then,
\begin{equation}
\label{pastroppres}
\Big\| 1 - |v| \Big\|_{L^\infty(\R^2)} \geq \frac{\varepsilon(v)^2}{10}.
\end{equation}
Moreover, there exists a universal constant $K_2 > 0$ such that
$$\varepsilon(v) \leq K_2 E(v).$$
\end{prop}

In particular, the solutions $u_\p$ given by Theorem \ref{dim2}, satisfy in view of Proposition \ref{firstprop},
$$\varepsilon_\p \equiv \varepsilon(u_\p) \to 0, \ {\rm as} \ \p \to 0,$$
so that we deal with a transonic limit. In \cite{IordSmi1,JoneRob1,JonPuRo1}, it is proposed to study this transonic limit of solutions $v$ in the new anisotropic space scale,
$$\tilde{x}_1 = \varepsilon(v) x_1, \ {\rm and} \ \tilde{x}_2 = \frac{\varepsilon(v)^2}{\sqrt{2}} x_2.$$
Considering the real-valued function
$$\eta \equiv 1 - |v|^2,$$
and performing the change of variables above, we introduce the rescaled map $N_v$ defined by
\begin{equation}
\label{scaling}
N_v(x) = \frac{6}{\varepsilon(v)^2} \eta \Big( \frac{x_1}{\varepsilon(v)}, \frac{\sqrt{2} x_2}{\varepsilon(v)^2} \Big).
\end{equation}
Notice that the same long-wave anisotropic scaling is performed to derive the Kadomtsev-Petviashvili equation, for instance in the water-wave context (see e.g. \cite{AbloSeg1,LannSau1}). It is formally shown in \cite{IordSmi1, JoneRob1, JonPuRo1} that the renormalized amplitude $N_v$ of solutions to \eqref{TWc} converges, as the speed $c(v)$ converges to $\sqrt{2}$, i.e. as $\varepsilon(v) \to 0$, to solitary wave solutions to the two-dimensional Kadomtsev-Petviashvili equation \eqref{KP}, that is
\renewcommand{\theequation}{KP I}
\begin{equation}
\label{KP}
\partial_t \psi + \psi \partial_1 \psi + \partial_1^3 \psi - \partial_1^{- 1} (\partial_2^2 \psi) = 0.
\end{equation}
Our main goal in this paper is to provide a rigorous mathematical proof of that convergence for the branch of minimizing solutions presented in Theorem \ref{dim2}.

Solitary waves are localized solutions to \eqref{KP} of the form $\psi(x, t) = w(x_1 - \sigma t, x_2)$, where $w$ belongs to the energy space for \eqref{KP}, i.e. the space $Y(\R^2)$ defined as the closure of $\partial_1 \boC_c^\infty(\R^2)$ for the norm
$$\| \partial_1 f \|_{Y(\R^2)} \equiv \Big( \| \nabla f \|_{L^2(\R^2)}^2 + \| \partial_1^2 f \|_{L^2(\R^2)}^2 \Big)^\frac{1}{2}.$$
The parameter $\sigma \geq 0$ denotes the speed of the solitary wave. The equation of a solitary wave $w$ of speed $\sigma = 1$ is given by
\renewcommand{\theequation}{SW}
\begin{equation}
\label{SW}
\partial_1 w - w \partial_1 w - \partial_1^3 w + \partial_1^{- 1} (\partial_2^2 w) = 0.
\end{equation}
When $w \in Y(\R^2)$, the function $\partial_1^{- 1} \partial_2 w$ is well-defined (see \cite{deBoSau1}), so that \eqref{SW} makes sense.

In contrast with the Gross-Pitaevskii equation, the range of speeds is the full positive axis. In particular, there are no solitary waves of negative speed (see \cite{deBoSau1}). Given any $\sigma \geq 0$, a solitary wave $w_\sigma$ of speed $\sigma$ is deduced from a solution $w$ to \eqref{SW} by the scaling
\renewcommand{\theequation}{\arabic{equation}}
\setcounter{equation}{6}
\begin{equation}
\label{scalingsw}
w_\sigma(x_1, x_2) = \sigma w(\sqrt{\sigma} x_1, \sigma x_2).
\end{equation}
Solitary waves may be obtained in dimension two minimizing the Hamiltonian keeping the $L^2$-norm fixed (see \cite{deBoSau3,deBoSau1}). Like \eqref{GP}, equation \eqref{KP} is indeed hamiltonian, with Hamiltonian given by
$$E_{KP}(\psi) = \frac{1}{2} \int_{\R^2} (\partial_1 \psi)^2 + \frac{1}{2} \int_{\R^2} (\partial_1^{-1}(\partial_2 \psi))^2 - \frac{1}{6} \int_{\R^2} \psi^3,$$
and the $L^2$-norm of $\psi$ is conserved as well. Setting
$$S(N) = E_{KP}(N) + \frac{\sigma}{2} \int_{\R^2} N^2,$$
we term ground state, a solitary wave $N$ which minimizes the action $S$ among all non-constant solitary waves of speed $\sigma$ (see \cite{deBoSau2} for more details). In dimension two, a solitary wave is a ground state if and only if it minimizes the Hamiltonian $E_{KP}$ keeping the $L^2$-norm fixed (see \cite{deBoSau3}). The constant $\boS_{KP}$, which appears in Theorem \ref{dim2}, denotes the action $S(N)$ of the ground states $N$ of speed $\sigma = 1$.

Going back to the solutions $u_\p$ of Theorem \ref{dim2}, we may drop the invariance under translations of our problem, assuming without loss of generality, since $|u_\p(x)| \to 1$, as $|x| \to + \infty$ (see \cite{Graveja1}), that $\eta_\p \equiv 1 - |u_\p|^2$ achieves its maximum at the origin, i.e.
$$\big\| \eta_\p \big\|_{L^\infty(\R^2)} = \big| \eta_\p(0) \big|.$$
We next consider the map
$$N_\p \equiv N_{u_\p}.$$
Notice that the origin is a maximum point for $N_\p$, and that in view of \eqref{pastroppres}, we have 
\begin{equation}
\label{pastriv}
N_\p(0) \geq \frac{3}{5}.
\end{equation}
Our main result is

\begin{theorem}
\label{convGPKP}
There exists a subsequence $(\p_n)_{n \in \N}$, tending to $0$ as $n$ tends to $+ \infty$, and a ground state $N_0$ of \eqref{KP} such that
$$N_{\p_n} \to N_0 \ {\rm in} \ W^{k, q}(\R^2), \ {\rm as} \ n \to + \infty,$$
for any $k \in \N$ and any $1 < q \leq + \infty$.
\end{theorem}

\begin{remark}
There is a well-known explicit solitary wave solution to \eqref{KP} of speed $1$, namely the so-called "lump" solution, which may be written as
$$w_\ell(x_1, x_2) = 24 \frac{3 - x_1^2 + x_2^2}{(3 + x_1^2 + x_2^2)^2}.$$
It is conjectured that the "lump" solution is a ground state. It is also conjectured that the ground state is unique, up to the invariances of the problem. If this was the case, then the full family $(N_{\p})_{\p > 0}$ would converge to $w_\ell$, as $\p \to 0$.
\end{remark}

So far, we have only discussed properties of the modulus of $u_\p$. However, in our argument, the phase is central as well. More precisely, if $\p$ is sufficiently small, then $u_\p$ has no zero in view of \eqref{pastroppres}, and we may lift it as $u_\p = \varrho_\p \exp i \varphi_\p$. Setting
\begin{equation}
\label{theta}
\Theta_\p(x) = \frac{6 \sqrt{2}}{\varepsilon_\p} \varphi_\p \Big( \frac{x_1}{\varepsilon_\p}, \frac{\sqrt{2} x_2}{\varepsilon_\p^2} \Big),
\end{equation}
we prove

\begin{prop}
\label{bangkok}
Let $(\p_n)_{n \in \N}$ and $N_0$ be as in Theorem \ref{convGPKP}. Passing possibly to a further subsequence, we have
$$\partial_1 \Theta _{\p_n} \to N_0 \ {\rm in} \ W^{k, q}(\R^2), \ {\rm as} \ n \to + \infty,$$
for any $k \in \N$ and any $1 < q \leq + \infty$.
\end{prop}

\begin{remark}
Equation \eqref{KP} is a higher dimensional extension of the well-known Korteweg-de Vries equation \eqref{KdV}, which may be written as
\renewcommand{\theequation}{KdV}
\begin{equation}
\label{KdV}
\partial_t \psi + \psi \partial_1 \psi + \partial_1^3 \psi = 0.
\end{equation}
In dimension one, travelling wave solutions $v_c$ to \eqref{TWc} are related to the classical soliton of the Korteweg-de Vries equation as follows. Setting $\varepsilon = \sqrt{2 - c^2}$, we consider the rescaled function
$$N_\varepsilon(x) = \frac{6}{\varepsilon^2} \eta_c \Big( \frac{x}{\varepsilon} \Big),$$
where $\eta_c \equiv 1 - |v_c|^2$. An explicit integration of \eqref{TWc} in dimension one leads to
$$N_\varepsilon(x) = N(x) \equiv \frac{3}{ \ch^2 \big( \frac{x}{2} \big)},$$
where $N$ is the classical soliton to the Korteweg-de-Vries equation. Concerning the phase $\varphi_c$ of $v_c$, we consider the scale change
$$\Theta_\varepsilon(x) = \frac{6\sqrt{2}}{\varepsilon} \varphi_c \Big( \frac{x}{\varepsilon} \Big),$$
so that we obtain similarly
$$ \Theta_\varepsilon(x)' = \sqrt {1 - \frac{\varepsilon^2}{2}} \frac{N(x)}{1 - \frac{\varepsilon^2}{2} N(x)} \to N(x), \ {\rm as} \ \varepsilon \to 0.$$
\end{remark}

\begin{remark}
Let $u_c$ be a solution to \eqref{TWc} in dimension three, which may be written as $u_c = \varrho_c \exp i \varphi_c$, and denote
$$N_c(x) = \frac{6}{\varepsilon^2} \eta_c \Big( \frac{x_1}{\varepsilon}, \frac{\sqrt{2} x_2}{\varepsilon^2}, \frac{\sqrt{2} x_3}{\varepsilon^2} \Big), \ {\rm and} \ \Theta_c(x) = \frac{6 \sqrt{2}}{\varepsilon} \varphi_c \Big( \frac{x_1}{\varepsilon}, \frac{\sqrt{2} x_2}{\varepsilon^2}, \frac{\sqrt{2} x_3}{\varepsilon^2} \Big),$$
where $\eta_c \equiv 1 - \varrho_c^2$ and $\varepsilon = \sqrt{2 - c^2}$. Then, it is also formally shown in \cite{IordSmi1,JoneRob1,JonPuRo1} that the functions $N_c$ and $\partial_1 \Theta_c$ converge, as the parameter $\varepsilon$ converges to $0$, to a solitary wave solution $w$ to the three-dimensional Kadomtsev-Petviashvili equation \eqref{KP}, which writes
$$\partial_t \psi + \psi \partial_1 \psi + \partial_1^3 \psi - \partial_1^{- 1} (\partial_2^2 \psi + \partial_3^2 \psi) = 0.$$
In particular, the equation for the solitary wave $w$ is now written as
$$\partial_1 w - w \partial_1 w - \partial_1^3 w + \partial_1^{- 1} (\partial_2^2 w + \partial_3^2 w) = 0.$$
However, the existence of a transonic branch of solutions is still an open problem, at least on the mathematical level. This branch of solutions is conjectured in \cite{IordSmi1,JoneRob1} in view of numerical computations and formal arguments.
\end{remark}

\subsection{Some elements in the proofs}

The first element in the proofs of Theorem \ref{convGPKP} and Proposition \ref{bangkok} deals with the asymptotic behaviour of $\varepsilon_\p$ as a function of $\p$.

\begin{lemma}[\cite{BetGrSa1}]
\label{ppetit} 
Let $\varepsilon_\p = \varepsilon(u_\p) = \sqrt{2 - c(u_\p)^2}$. There exist some positive constants $K_3$ and $K_4$, not depending on $\p$, such that
\renewcommand{\theequation}{\arabic{equation}}
\setcounter{equation}{9}
\begin{equation}
\label{epsilonetp}
K_3 \p \leq \varepsilon_\p \leq K_4 \p,
\end{equation}
for any $\p$ sufficiently small.
\end{lemma}

The second step is to derive estimates on the renormalized maps $N_\p$, which do not depend on $\p$. More precisely, we prove

\begin{prop}
\label{Sobolevbound}
Let $k \in \N$ and $1 < q \leq + \infty$. There exists some constant $K(k, q)$, depending possibly on $k$ and $q$, but not on $\p$, such that
\begin{equation}
\label{higher}
\| N_\p \|_{W^{k, q}(\R^2)} +\| \partial_1 \Theta_\p \|_{W^{k, q}(\R^2)}+\eps_\p\| \partial_2 \Theta_\p \|_{W^{k, q}(\R^2)} \leq K(k, q),
\end{equation}
for any $\p$ sufficiently small.
\end{prop}

At this stage, we may invoke standard compactness theorems to assert that there exists some subsequence $(\p_n)_{n \in \N}$, tending to $0$ as $n$ tends to $+ \infty$, and a function $N_0$ such that, for any $k \in \N$ and any compact subset $K$ of $\R^2$,
$$N_{\p_n} \to N_0 \ {\rm in} \ \boC^k(K), \ {\rm as} \ n \to + \infty.$$
In view of \eqref{pastriv}, we have
$$N_0(0) \geq \frac{3}{5},$$
so that $N_0$ is not identically constant. Moreover, we also have

\begin{lemma}
\label{roti}
The function $N_0$ is a non-constant solution to \eqref{SW}.
\end{lemma}

In order to complete the proof of Theorem \ref{convGPKP}, it remains to establish strong convergence on the whole plane. For this last step, we essentially rely on a variational argument, proving a kind of gamma-convergence of the energies, combined with a concentration-compactness result for constrained minimizers of \eqref{KP} established in \cite{deBoSau3}.

As a matter of fact, considering scalings \eqref{scaling} and \eqref{theta}, the momentum $p(u_\p)$ can be expressed as
$$p(u_\p) = \frac{\varepsilon_\p}{72} \int_{\R^2} N_\p \partial_1 \Theta_\p,$$
while the energy $E(u_\p)$ has the expansion
$$E(u_\p) = \sqrt{2} \frac{\varepsilon_\p}{144} \Big( E_0(N_\p, \Theta_\p) + \varepsilon_\p^2 E_2(N_\p, \Theta_\p) + \varepsilon_\p^4 E_ 4 (N_\p, \Theta_\p) \Big).$$
It turns out that the functions $E_0$, $E_2$ and $E_4$ are uniformly bounded for $\p$ approaching $0$. Moreover, $E_0$ and $E_2$ are given by the expressions
$$E_0(N_\p, \Theta_\p) = \int_{\R^2} \Big( N_\p^2 + (\partial_1 \Theta_\p)^2 \Big),$$
and
\begin{equation}
\label{defE20}
E_2(N_\p, \Theta_\p) = \int_{\R^2} \Big( \frac{1}{2} (\partial_1 N_\p)^2 + \frac{1}{2} (\partial_2 \Theta_\p)^2 - \frac{1}{6} N_\p (\partial_1 \Theta_\p)^2 \Big). 
\end{equation}
In the course of our proof, we will show that
\begin{equation}
\label{sim}
N_\p \sim \partial_1 \Theta_\p, \ {\rm as} \ \p \to 0,
\end{equation}
and that the difference is actually of order $\varepsilon_\p^2$. This yields, at least heuristically,
$$p(u_\p)\sim \frac{\varepsilon_\p}{72}\int_{\R^2} N_\p^2, \ {\rm and} \ E(u_\p) \sim 
\sqrt{2} \frac{\varepsilon_\p}{72}\int_{\R^2} N_\p^2\sim \sqrt{2}p(u_\p),$$
so that the discrepancy term
$$\Sigma(u_\p) = \sqrt{2} p(u_\p) - E(u_\p),$$
tends to $0$ as $\p \to + \infty$.

The \eqref{KP} energy appears when we consider the second order term. Inserting at least formally relation \eqref{sim} into \eqref{defE20}, we are led to
\begin{equation}
\label{equivalence}
E_2(N_\p, \Theta_\p) \sim E_{KP}(N_\p), \ {\rm as} \ \p \to 0.
\end{equation}
Using some precise estimates on the solutions, we will actually show that
\begin{equation}
\label{equivalence2}
E_2(N_\p, \Theta_\p) \sim E_{KP}(\partial_1\Theta_\p), \ {\rm as} \ \p \to 0,
\end{equation}
since it turns out that it is easier to work, in view of the nonlocal term in the \eqref{KP} energy, with $\partial_1 \Theta_\p$ than with $N_\p$, these two terms having the same limit in view of \eqref{sim}.

The proof of \eqref{equivalence2} amounts to a careful analysis of any lower order terms, including terms provided by $E_0$. In particular, we obtain for the discrepancy functional,

\begin{lemma}
\label{discrepancy}
We have
\begin{equation}
\label{sigma}
\Sigma(u_\p) = - \frac{\sqrt{2} \varepsilon_\p^3}{144} E_{KP}(\partial_1 \Theta_\p) + \underset{\p \to 0}{o} \big( \varepsilon_\p^3 \big).
\end{equation}
\end{lemma}

We then use the lower bound on $\Sigma(u_\p)$ provided by the left-hand side of \eqref{estimE2} to derive a precise upper bound on $E_{KP}(\partial_1 \Theta_\p)$. More precisely, we show

\begin{lemma}
\label{upperbound}
We have
\begin{equation}
\label{page}
- \frac{1}{54 \boS_{KP}^2} \bigg( \int_{\R^2} \big( \partial_1 \Theta_\p \big)^2 \bigg)^3 \leq E_{KP}(\partial_1 \Theta_\p) \leq - \frac{1}{54 \boS_{KP}^2} \bigg( \int_{\R^2} \big( \partial_1 \Theta_\p \big)^2 \bigg)^3 + \underset{\p \to 0}{o}(1).
\end{equation}
\end{lemma}

In particular, the function $\partial_1 \Theta_\p$, or alternatively $N_\p$, has approximatively the energy of a ground state for \eqref{KP} corresponding to its $L^2$-norm. The proof of Theorem \ref{convGPKP} is then completed using a concentration-compactness argument of \cite{deBoSau3}. This result yields the strong convergence of some subsequence $(\partial_1 \Theta_{\p_n})_{n \in \N}$ in the space $Y(\R^2)$.

\begin{prop}
\label{Yconv}
There exists a subsequence $(\p_n)_{n \in \N}$, tending to $0$ as $n$ tends to $+ \infty$, and a ground state $N_0$ of \eqref{KP} such that
$$\partial_1 \Theta_{\p_n} \to N_0 \ {\rm in} \ Y(\R^2), \ {\rm and} \ N_{\p_n} \to N_0 \ {\rm in} \ L^2(\R^2), \ {\rm as} \ n \to + \infty.$$
\end{prop}

In order to improve the convergence, we finally invoke the estimates of Proposition \ref{Sobolevbound}. This concludes the proofs of Theorem \ref{convGPKP} and Proposition \ref{bangkok} giving the convergence in any space $W^{k, q}(\R^2)$ by standard interpolation theory.

To conclude this introduction, let us emphasize that the results in this paper only concern travelling waves. This raises quite naturally the corresponding issue for the time-dependent equations. More precisely, in which sense do the Korteweg-de Vries equation in dimension one and the Kadomtsev-Petviashvili equation in higher dimensions approximate the Gross-Pitaevskii equation in the transonic limit ? Notice that this question has already been formally addressed in the one-dimensional case in \cite{KuznZak1}.

\subsection{Outline of the paper}

The paper is organized as follows. Sections \ref{PropKP} and \ref{PropGP} are devoted to various properties of solitary wave solutions to \eqref{KP} and travelling wave solutions to \eqref{TWc} which are subsequently used. In Section \ref{Lent}, we perform the expansion of \eqref{TWc} with respect to the small parameter $\varepsilon$ occurring in the definition of the slow space variables. Terms in this expansion are more clearly analyzed in Fourier variables. Various kernels then appear, which are studied in Section \ref{Linear}. In Section \ref{Unity}, we provide Sobolev bounds on $N_\p$ and prove Proposition \ref{Sobolevbound}. Finally, we prove our main theorems in Section \ref{Concentration}.

\numberwithin{cor}{section}
\numberwithin{equation}{section}
\numberwithin{lemma}{section}
\numberwithin{prop}{section}
\numberwithin{remark}{section}
\numberwithin{theorem}{section}
\section{Some properties of solitary wave solutions to \eqref{KP}}
\label{PropKP}

We first recall some facts about equation \eqref{KP}, which will enter in some places in our proofs.

\subsection{Rewriting the solitary wave equation}
\label{newform}

The existence and qualitative properties of the solutions $w$ to \eqref{SW} in the energy space $Y(\R^2)$ are considered in the series of papers \cite{deBoSau1,deBoSau2,deBoSau3}. In \cite{deBoSau2}, a new formulation of \eqref{SW} is provided which turns out to be also fruitful in our context. Applying the operator $\partial_1$ to \eqref{SW}, we obtain
\begin{equation}
\label{eqSW2}
\partial_1^4 w - \Delta w + \frac{1}{2} \partial_1^2 (w^2) = 0.
\end{equation}
The Fourier transform of \eqref{eqSW2} has the following simple form
\begin{equation}
\label{fouSW}
\widehat{w}(\xi) = \frac{1}{2} \frac{\xi_1^2}{|\xi|^2 + \xi_1^4} \widehat{w^2}(\xi),
\end{equation}
so that we may recast \eqref{eqSW2} as a convolution equation
\begin{equation}
\label{convSW}
w = \frac{1}{2} K_0 \star w^2,
\end{equation}
where the Fourier transform of the kernel $K_0$ is given by
\begin{equation}
\label{kernel0}
\widehat{K_0}(\xi) = \frac{\xi_1^2}{|\xi|^2 + \xi_1^4}.
\end{equation}
In view of \eqref{fouSW}, equation \eqref{convSW} provides an equivalent formulation to \eqref{SW}, i.e. any solution $w$ to \eqref{convSW} in the energy space $Y(\R^2)$ is also solution to \eqref{SW}.

Several properties of the kernel $K_0$ are studied in \cite{Graveja7}. In particular, it is proved there that $K_0$ belongs to $L^p(\R^2)$ for any $1 < p < 3$ (see also Lemma \ref{noyaudur}).

\subsection{Existence of ground state solutions}
\label{GroundKP}

Given any $\mu \geq 0$, the minimization problem
\renewcommand{\theequation}{$\boP_{KP}(\mu)$}
\begin{equation}
\label{minikp}
\boE^{KP}_{\min}(\mu) = \inf \Big\{ E_{KP}(w), w \in Y(\R^2), \int_{\R^2} |w|^2 = \mu \Big\},
\end{equation}
is considered in \cite{deBoSau3}, where the existence of minimizers is established. The minimizers $N$ for this problem happen to be ground states for \eqref{KP}. They are solutions to
\renewcommand{\theequation}{\arabic{equation}}
\numberwithin{equation}{section}
\setcounter{equation}{4}
\begin{equation}
\label{SWsigma}
\sigma \partial_1 N - N \partial_1 N - \partial_1^3 N + \partial_1^{- 1} (\partial_2^2 N) = 0.
\end{equation}
The speed $\sigma$ appears as a Lagrange multiplier associated to \eqref{minikp}. In particular, $\sigma$ is not necessarily equal to $1$. The proof in \cite{deBoSau3} relies on the following concentration-compactness result, which gives the compactness of minimizing sequences to \eqref{minikp}.

\begin{theorem}[\cite{deBoSau3}]
\label{ConcKP}
Let $\mu \geq 0$, and let $(w_n)_{n \in \N}$ be a minimizing sequence to \eqref{minikp} in $Y(\R^2)$. Then, there exist some points $(a_{n})_{n \in \N}$ and a function $N \in Y(\R^2)$ such that, up to some subsequence,
$$w_n(\cdot - a_n) \to N \ {\rm in} \ Y(\R^2), \ {\rm as} \ n \to + \infty.$$
The limit function $N$ is solution to the minimization problem \eqref{minikp}. In particular, $N$ is a ground state for \eqref{KP}.
\end{theorem}

\subsection{Scale invariance}
\label{armstrong}

As mentioned in the introduction, if $w$ is solution to \eqref{SW}, then, for any $\sigma > 0$, the map $w_\sigma$ defined by \eqref{scalingsw} is solution to \eqref{SWsigma}, i.e. $w_\sigma$ is a solitary wave solution to \eqref{KP} with speed $\sigma$. Concerning the energy, we notice that
$$\int_{\R^2} |w_\sigma|^2 = \sqrt{\sigma} \int_{\R^2} |w|^2, \ \int_{\R^2} |w_\sigma|^3 = \sigma^\frac{3}{2} \int_{\R^2} |w|^3, \ \int_{\R^2} |\partial_1 w_\sigma|^2 = \sigma^\frac{3}{2} \int_{\R^2} |\partial_1 w|^2,$$
and
$$\int_{\R^2} \Big( \partial_1^{-1}(\partial_2 w_\sigma) \Big)^2 = \sigma^\frac{3}{2} \int_{\R^2} \Big( \partial_1^{-1}(\partial_2 w) \Big)^2.$$
It follows that
\begin{equation}
\label{scaledenergy}
E_{KP}(w_\sigma) = \sigma^\frac{3}{2} E_{KP}(w), \ {\rm and} \ \int_{\R^2} |w_\sigma|^2 = \sqrt{\sigma} \int_{\R^2} |w|^2.
\end{equation}
It is shown in \cite{deBoSau3} that ground states $N$ with speed $\sigma = 1$ correspond to solutions to \eqref{minikp} for
$$\mu = \mu^* \equiv 3 \boS_{KP}.$$
As a matter of fact, it is proved in \cite{deBoSau1,Graveja7} that any solution $w$ to \eqref{SW} satisfies the relations 
$$E_{KP}(w) = - \frac{1}{6} \int_{\R^2} w^2, \ {\rm and} \ S(w) = \frac{1}{3} \int_{\R^2} w^2,$$
so that the energy and the $L^2$-norm of ground states $N$ with speed $\sigma = 1$ are given by
$$E_{KP}(N) = - \frac{1}{2} \boS_{KP}, \ {\rm and} \ \int_{\R^2} N^2 = 3 \boS_{KP} = \mu^*.$$
Relations \eqref{scaledenergy} then provide

\begin{lemma}
\label{MiniKP}
Let $N \in Y(\R^2)$. Given any $\sigma \geq 0$, the map $N_\sigma$ defined by \eqref{scalingsw} is a minimizer for $\boE^{KP}_{\min}(\sqrt{\sigma} \mu^*)$ if and only if $N$ is a minimizer for $\boE^{KP}_{\min}(\mu^*)$. In particular, we have
\begin{equation}
\label{minKP}
\boE^{KP}_{\min}(\mu) = - \frac{\mu^3}{54 \boS_{KP}^2}, \ \forall \mu \geq 0.
\end{equation}
Moreover, $N_\sigma$ and $N$ are ground states for \eqref{KP}, with speed $\sigma$, respectively, $1$. In particular, they are solutions to \eqref{SWsigma}, respectively, \eqref{SW}.
\end{lemma}

\begin{proof}
Given any $\mu > 0$, we denote $\Lambda^2_\mu(\R^2) = \{ w \in L^2(\R^2), \ {\rm s.t.} \ \int_{\R^2} |w|^2 = \mu \}$. In view of \eqref{scaledenergy}, the function $w \mapsto w_\sigma$ maps $\Lambda^2_{\mu^*}(\R^2)$ onto $\Lambda^2_{\mu^* \sqrt{\sigma}}(\R^2)$, such that
$$E_{KP}(w_\sigma) = \sigma^\frac{3}{2} E_{KP}(w).$$
Hence, $N_\sigma$ is a minimizer for $\boE^{KP}_{\min}(\mu^* \sqrt{\sigma})$ if and only if $N$ is a minimizer for $\boE^{KP}_{\min}(\mu^*)$. Moreover,
$$\boE^{KP}_{\min}(\mu^* \sqrt{\sigma}) = \sigma^\frac{3}{2} \boE^{KP}_{\min}(\mu^*) = -\frac{\sigma^\frac{3}{2} \boS_{KP}}{2}.$$
Identity \eqref{minKP} follows letting $\sigma = \frac{\mu^2}{(\mu^*)^2}$. The last statements of Lemma \ref{MiniKP} are proved in \cite{deBoSau3}.
\end{proof}

In the course of our proofs, we will encounter sequences $(w_n)_{n \in \N}$ which are not exactly minimizing sequences for \eqref{minikp}, but which satisfy
\begin{equation}
\label{miniseq}
E_{KP}(w_n) \to \boE^{KP}_{\min}(\mu), \ {\rm and} \ \int_{\R^2} w_n^2 \to \mu, \ {\rm as} \ n \to + \infty,
\end{equation}
for some positive number $\mu$. In this case, we will invoke the following variant (and in fact, consequence) of Theorem \ref{ConcKP}.

\begin{prop}
\label{Goodconc}
Let $\mu_0 > 0$, and $(w_n)_{n \in \N}$ denote a sequence of functions in $Y(\R^2)$ satisfying \eqref{miniseq} for $\mu = \mu_0$. Then, there exist some points $(a_n)_{n \in \N}$ and a ground state solution $N_\sigma$ to \eqref{SWsigma}, with $\sigma = \frac{\mu_0^2}{(\mu^*)^2}$, such that, up to some subsequence,
$$w_n(\cdot - a_n) \to N_\sigma \ {\rm in} \ Y(\R^2), \ {\rm as} \ n \to + \infty.$$
\end{prop}

\begin{proof}
We denote
$$\mu_n = \int_{\R^2} w_n^2, \ {\rm and} \ \sigma_n = \frac{\mu_0^2}{\mu_n^2},$$
and consider the functions
$$z_n(x_1, x_2) = \sigma_n w_n(\sqrt{\sigma_n} x_1, \sigma_n x_2).$$
In view of \eqref{scaledenergy} and \eqref{miniseq},
\begin{equation}
\label{rodeo}
\sigma_n \to 1, \ {\rm as} \ n \to + \infty,
\end{equation}
and $(z_n)_{n \in \N}$ is a minimizing sequence of \eqref{minikp} for $\mu = \mu_0$. Therefore, by Theorem \ref{ConcKP}, there exist some points $(a_n)_{n \in \N}$ and a minimizer $N_\sigma$ to \eqref{minikp} for $\mu = \mu_0$ such that, up to some subsequence,
\begin{equation}
\label{lasso}
z_n(\cdot - a_n) \to N_\sigma \ {\rm in} \ Y(\R^2), \ {\rm as} \ n \to + \infty.
\end{equation}
In particular, it follows from Lemma \ref{MiniKP} that $N_\sigma$ is solution to \eqref{SWsigma}, with $\sigma = \frac{\mu_0^2}{(\mu^*)^2}$.
We now denote
$$N_n(x_1, x_2) = \frac{1}{\sigma_n} N_\sigma \Big( \frac{x_1}{\sqrt{\sigma_n}}, \frac{x_2}{\sigma_n} \Big),$$
so that, by the change of variables $(y_1, y_2) = (\sqrt{\sigma_n} x_1, \sigma_n x_2)$,
\begin{align*}
\| z_n(\cdot - a_n) - N_\sigma \|_{Y(\R^2)}^2 = & \sqrt{\sigma_n} \| w_n(\cdot - a_n) - N_n \|_{L^2(\R^2)}^2 + \sigma_n^\frac{3}{2} \| \partial_1 w_n(\cdot - a_n) - \partial_1 N_n \|_{L^2(\R^2)}^2\\
+ & \sigma_n^\frac{3}{2} \| \partial_1^{-1} \partial_2 w_n(\cdot - a_n) - \partial_1^{-1} \partial_2 N_n \|_{L^2(\R^2)}^2.
\end{align*}
By \eqref{rodeo} and \eqref{lasso}, we have
$$w_n(\cdot - a_n) - N_n \to 0 \ {\rm in} \ Y(\R^2), \ {\rm as} \ n \ \to + \infty.$$
Proposition \ref{Goodconc} follows provided we first prove that
$$N_n \to N_\sigma \ {\rm in} \ Y(\R^2), \ {\rm as} \ n \ \to + \infty.$$
This last assertion is itself a consequence of the general observation that
$$\lambda \psi \big( \sqrt{\mu} \cdot, \mu \cdot \big) \to \psi \ {\rm in} \ L^2(\R^2), \ {\rm as} \ \lambda \to 1 \ {\rm and} \ \mu \to 1,$$
which may be deduced from the dominated convergence theorem, when $\psi$ is in $\boC_c^\infty(\R^2)$, then, using the density of $\boC_c^\infty(\R^2)$ into $L^2(\R^2)$, when $\psi$ only belongs to $L^2(\R^2)$.
\end{proof}

\section{Some properties of solutions to \eqref{TWc}}
\label{PropGP}

In this section, we gather a number of properties of solutions to \eqref{TWc}, which enter in our asymptotic analysis. Most of these results are available in the literature on the subject. 

\subsection{General solutions}

Let $v$ be a finite energy solution to \eqref{TWc} on $\R^2$. It can be shown using various elliptic estimates (see \cite{Farina1,Tarquin1,BetGrSa1}) that there exists some positive constant $K$, not depending on $c$, such that
\begin{equation}
\label {elinfiniv}
\Big\| 1 - |v| \Big\|_{L^\infty(\R^2)} \leq 1,
\end{equation}
and
\begin{equation}
\label{elinfinigrad}
\| \nabla v \|_{L^\infty(\R^2)} \leq K \Big( 1 + \frac{c^2}{4} \Big)^\frac{3}{2}.
\end{equation}
In view of \eqref{subsonicwave}, estimates \eqref{elinfiniv} and \eqref{elinfinigrad} may be recast as
\begin{equation}
\label{elinfiniestimates}
\| \eta \|_{L^\infty(\R^2)} + \| \nabla v \|_{L^\infty(\R^2)} \leq K,
\end{equation}
where we have set $\eta \equiv 1 - |v|^2$. For higher order derivatives, it similarly follows from the proof of Lemma 2.1 in \cite{BetGrSa1} that there exists some positive constant $K(k)$, not depending on $c$, such that
\begin{equation}
\label{ckestimates}
\| v \|_{\boC^k(\R^2)} \leq K(k),
\end{equation}
for any $k \in \N$.

More generally, we have
\begin{equation}
\label{sobestimates}
\| \eta \|_{W^{k, q}(\R^2)} + \| \nabla v \|_{W^{k, q}(\R^2)} \leq K(c, k, q),
\end{equation}
for any $k \in \N$ and any $1 < q < + \infty$ (see \cite{Graveja3}). Notice that the constant $K(c, k, q)$ possibly depends on the speed $c$, so that we may have
$$K(c, k, q) \to + \infty, \ {\rm as} \ c \to \sqrt{2}.$$
Before establishing the convergence of the rescaled functions $N_\p$ and $\Theta_\p$, we shall need to establish their boundedness in the spaces $W^{k, q}(\R^2)$. This requires to get some control upon the dependence on $c$ of the constant $K(c, k, q)$. The proof of Proposition \ref{Sobolevbound} in Section \ref{Unity} below provides such a control.

We will also take advantage of the fact that the maps $u_\p$ have small energy. Indeed, in view of \eqref{subsonicwave} and elliptic estimate \eqref{elinfiniestimates}, we may show that, if a solution $v$ to \eqref{TWc} has sufficiently small energy, it does not vanish. More precisely, we have

\begin{lemma}[\cite{BetGrSa1}]
\label{nozero}
There exists a universal constant $E_0$ such that, if $v$ is a solution to \eqref{TWc} which satisfies $E(v) \leq E_0$, then
\begin{equation}
\label{loindezero}
\frac{1}{2} \leq |v| \leq 2.
\end{equation}
\end{lemma}

If $v$ satisfies \eqref{loindezero}, then we may lift it as
$$v = \varrho \exp i \varphi,$$
where $\varphi$ is a real-valued, smooth function on $\R^2$ defined modulo a multiple of $2 \pi$. We have in that case,
$$\partial_j v = \Big( i \varrho \partial_j \varphi + \partial_j \varrho \Big) \exp i \varphi,$$
so that
\begin{equation}
\label{relevenergy0}
\langle i \partial_1 v, v \rangle = - \varrho^2 \partial_1 \varphi, \ {\rm and} \ e(v) = \frac{1}{2} \Big( |\nabla \varrho|^2 + \varrho^2 |\nabla \varphi|^2 \Big) + \frac{1}{4} \eta ^2.
\end{equation}
Moreover, the momentum $p$ takes the simple form
$$p(v) = \frac{1}{2} \int_{\R^2} \eta \partial_1 \varphi.$$
The system of equations for $\varrho$ and $\varphi$ is written as
\begin{equation}
\label{PolTWc1}
\frac{c}{2} \partial_1 \varrho^2 + \div \Big( \varrho^2 \nabla \varphi \Big) = 0,
\end{equation}
and
\begin{equation}
\label{PolTWc2}
c \varrho \partial_1 \varphi - \Delta \varrho - \varrho \Big( 1 - \varrho^2 \Big) + \varrho |\nabla \varphi|^2 = 0.
\end{equation}
Combining both the equations, the quantity $\eta$ satisfies
$$\Delta^2 \eta - 2 \Delta \eta + c^2 \partial_1^2 \eta= - 2 \Delta \big( |\nabla v|^2 + \eta^2 - c \eta \partial_1 \varphi \big) - 2 c \partial_1 \div \big( \eta \nabla \varphi \big),$$
where the left-hand side is linear with respect to $\eta$, whereas the right-hand side is (almost) quadratic with respect to $\eta$ and $\nabla \varphi$.

Multiplying \eqref{PolTWc1} by $\varphi$ and integrating by parts, we obtain a first relation for the momentum
\begin{equation}
\label{merida}
c p(v) = \int_{\R^N} \varrho^2 |\nabla \varphi|^2.
\end{equation}
In another direction, Pohozaev identities yield
\begin{equation}
\label{povrose1}
E(v) = \int_{\R^2} |\partial_1 v|^2, \ {\rm and} \ E(v) = \int_{\R^2} |\partial_2 v|^2 + c p(v).
\end{equation}
Introducing the quantities $\Sigma(v) = \sqrt{2} p(v)- E(v)$, the second identity in \eqref{povrose1} may be recast as
\begin{equation}
\label{acapulco}
\int_{\R^2} |\partial_2 v|^2 + \Sigma(v) = \Big( \sqrt{2} - \sqrt{2-\varepsilon(v)^2} \Big) p(v) = \frac{\varepsilon(v)^2}{\sqrt{2} + \sqrt{2 - \varepsilon(v)^2}} p(v).
\end{equation}
In the case $\Sigma(v) > 0$, this yields an interesting estimate for the transversal derivative $\partial_2 v$. Adding both the equalities in \eqref{povrose1}, we also derive a second relation for the momentum
$$\frac{1}{2} \int_{\R^2}\eta^2 = c p(v).$$
With similar arguments and combining with \eqref{merida}, we are led to

\begin{lemma}[\cite{BetGrSa1}]
\label{tropbien}
Let $v$ be a finite energy solution to \eqref{TWc} on $\R^2$ satisfying \eqref{loindezero}. Then, we have the identities
\begin{equation}
\label{tropbien1}
\Sigma(v) + \frac{1}{2} \int_{\R^2} |\nabla \varrho|^2 = \frac{\varepsilon(v)^2}{\sqrt{2} + c(v)} p(v),
\end{equation}
\begin{equation}
\label{tropbien2}
\int_{\R^2} |\nabla \varrho|^2 \Big( 1 + \frac{1}{\varrho^2} \Big) = \int_{\R^2} \eta |\nabla \varphi|^2,
\end{equation}
and the inequality
\begin{equation}
\label{fastoche}
E(v)\leq 7 c(v)^2 \int_{\R^2} \eta^2.
\end{equation}
\end{lemma}

In view of definition \eqref{GLE}, we have
$$\int_{\R^2} \eta^2 \leq 4 E(v),$$
so that inequality \eqref{fastoche} shows that the energy is comparable to the integral of $\eta^2$ for any solutions $v$ satisfying \eqref{loindezero}. When $\Sigma(v) > 0$, identity \eqref{tropbien1} shows that
$$\Sigma(v) \leq \frac{\varepsilon(v)^2}{\sqrt{2}} p(v) \leq K E(v)^2 p(v) \leq 2 K p(v)^3,$$
where we have invoked Proposition \ref{firstprop} for the second inequality. In particular, we obtain
$$E(v) \sim \sqrt{2}p(v),$$
as $E(v)$, or $p(v)$, approaches $0$.

In several places (in particular, in the proof of Proposition \ref{Sobolevbound}), we shall need estimates for higher order derivatives. For that purpose, we shall use

\begin{lemma}
\label{Derivphi}
Let $1 < q < + \infty$, and let $v$ be a finite energy solution to \eqref{TWc} on $\R^2$ satisfying \eqref{loindezero}. Then, there exists some constant $K(q)$, not depending on $c$, such that
\begin{equation}
\label{pepe}
\| \nabla \varphi \|_{L^q(\R^2)} \leq K(q) \| \eta \|_{L^q(\R^2)},\\
\end{equation}
More generally, given any index $\alpha = (\alpha_1, \alpha_2) \in \N^2$, there exist some constants $K(q, \alpha)$, not depending on $c$, such that
\begin{equation}
\label{meme}
\| \partial^\alpha (\nabla \varphi) \|_{L^q(\R^2)} \leq K(q, \alpha) \Big( \big\| \partial^\alpha \eta \big\|_{L^q(\R^2)} + \sum_{0 \leq \beta < \alpha} \| \partial^\beta \eta \|_{L^\infty(\R^2)} \| \partial^{\alpha - \beta} \big( \nabla \varphi \big) \|_{L^q(\R^2)} \Big).
\end{equation}
\end{lemma}

\begin{proof}
First notice that in view of \eqref{ckestimates} and \eqref{sobestimates}, the functions $\eta$ and $\nabla \varphi$ belong to $W^{k, q}(\R^2)$ for any $k \in \N$ and any $1 < q \leq + \infty$. In particular, the norms in inequalities \eqref{pepe} and \eqref{meme} are well-defined and finite. Lemma \ref{Derivphi} is then a consequence of the elliptic nature of equation \eqref{PolTWc1}, which may be written as
$$\Delta \varphi = \frac{c}{2} \partial_1 \eta + \div \big( \eta \nabla \varphi \big),$$
so that, more generally,
\begin{equation}
\label{alphaphi}
\Delta (\partial^\alpha \varphi) = \frac{c}{2} \partial_1 \partial^\alpha \eta + \div \big( \partial^\alpha (\eta \nabla \varphi) \big),
\end{equation}
for any $\alpha \in \N^2$.
Using standard elliptic estimates and inequality \eqref{subsonicwave}, we derive from \eqref{alphaphi} that
\begin{equation}
\label{elliptic}
\| \nabla (\partial^\alpha \varphi) \|_{L^q(\R^2)} \leq K(q) \Big( \big\| \partial^\alpha \eta \big\|_{L^q(\R^2)} + \big\| \partial^\alpha \big( \eta \nabla \varphi \big) \big\|_{L^q(\R^2)} \Big).
\end{equation}
For $\alpha = (0, 0)$, inequality \eqref{pepe} is a direct consequence of \eqref{elliptic} invoking \eqref{elinfiniestimates}. For $\alpha \neq (0, 0)$, the derivative $\partial^\alpha (\eta \nabla \varphi)$ may be written as
$$\partial^\alpha (\eta \nabla \varphi) = \sum_{0 \leq \beta \leq \alpha} \binom{\alpha}{\beta} \partial^\beta \eta \partial^{\alpha - \beta} \big( \nabla \varphi \big),$$
by Leibniz formula, so that
$$\| \partial^\alpha (\eta \nabla \varphi) \|_{L^q(\R^2)} \leq K(q, \alpha) \bigg( \| \partial^\alpha \eta \|_{L^q(\R^2)} \| \nabla \varphi \|_{L^\infty(\R^2)} + \sum_{0 \leq \beta < \alpha} \| \partial^\beta \eta \|_{L^\infty(\R^2)} \| \partial^{\alpha - \beta} \big( \nabla \varphi \big) \|_{L^q(\R^2)} \bigg).$$
Estimate \eqref{meme} follows from \eqref{elliptic} using again uniform bound \eqref{elinfiniestimates}.
\end{proof}

\subsection{Properties of $u_\p$}

We now restrict ourselves to the solutions $u_\p$ provided by Theorem \ref{dim2}. We begin with the

\begin{proof}[Proof of Lemma \ref{ppetit}]
In view of \eqref{estimE2}, we have
$$\Sigma_\p \equiv \Sigma(u_\p) \geq \frac{48 \sqrt{2}}{\boS_{KP}^2} \p^3 - K_0 \p^4,$$
for any $\p$ sufficiently small, whereas, by \eqref{tropbien1},
$$\Sigma_\p \leq \p \frac{\varepsilon_{\p}^2}{\sqrt{2}},$$
so that, combining both the inequalities, we obtain
$$\varepsilon_\p \geq \frac{9}{\boS_{KP}} \p.$$
On the other hand, in view of Proposition \ref{firstprop}, we have
$$\varepsilon_\p \leq K E_\p,$$
where we have set $E_\p \equiv E(u_\p)$. Since $E_\p \leq \sqrt{2} \p$, we conclude that \eqref{epsilonetp} holds. Moreover, we also have
\begin{equation}
\label{euhh}
K_5 E_\p \leq \varepsilon_{\p} \leq K_6 E_\p,
\end{equation}
for any $\p$ sufficiently small, and some positive constants $K_5$ and $K_6$, not depending on $\p$.
\end{proof}

Finally, since $\Sigma_\p > 0$ by \eqref{estimE2}, we deduce from Lemma \ref{ppetit} that \eqref{acapulco}, \eqref{tropbien1} and \eqref{tropbien2} may be recast as
\begin{equation}
\label{deriveta}
\int_{\R^2} \Big( |\nabla \varrho_\p|^2 + (\partial_2 u_\p)^2 \Big) + \bigg| \int_{\R^2} \eta_\p |\nabla \varphi_\p|^2 \bigg| \leq K \p^3,
\end{equation}
where we denote $u_\p = \varrho_\p \exp i \varphi_\p$. Since $(\partial_2 u_\p)^2 = \varrho_\p^2 (\partial_2 \varphi_\p)^2 + ( \partial_2 \varrho_\p)^2$ and $|\eta_\p| \leq 3 \varrho_\p^2$, we deduce that
$$\int_{\R^2} |\eta_\p| (\partial_2 \varphi_\p)^2 \leq 3 \int_{\R^2} \varrho_\p^2 (\partial_2 \varphi_\p)^2 \leq K \p^3,$$
so that
\begin{equation}
\label{deriveta2}
\bigg| \int_{\R^2} \eta_\p (\partial_1 \varphi_\p)^2 \bigg| \leq K \p^3.
\end{equation}

\section{\eqref{TWc} in the slow space variables}
\label{Lent}

\subsection{Expansion of the energy functionals}

In this subsection, we consider a finite energy map $v$ on $\R^2$, satisfying \eqref{loindezero}, and a small given parameter $\varepsilon > 0$. In view of assumption \eqref{loindezero}, we may lift $v$ as $v = \varrho \exp i \varphi$. Following the expansion given in the physical literature, we introduce anisotropic slow space variables $\tilde{x}_1 = \varepsilon x_1$, and $\tilde{x}_2 = \frac{\varepsilon^2}{\sqrt{2}} x_2$. We then consider the rescaled functions $N = N_{v, \varepsilon}$ and $\Theta = \Theta_{v, \varepsilon}$ defined as follows 
\begin{equation}
\label{slow-var}
N(x) = \frac{6}{\varepsilon^2} \eta \Big( \frac{x_1}{\varepsilon}, \frac{\sqrt{2} x_2}{\varepsilon^2} \Big), \ {\rm and} \ \Theta(x) = \frac{6 \sqrt{2}}{\varepsilon} \varphi \Big( \frac{x_1}{\varepsilon}, \frac{\sqrt{2} x_2}{\varepsilon^2} \Big).
\end{equation}
We next express the functionals $p$ and $E$ in terms of the functions $N$, $\Theta$ and $\varepsilon$. In the course of the analysis, we will also compute several other integral quantities in the rescaled variables. For instance,
$$\int_{\R^2} N^2 = \frac{18 \sqrt{2}}{\varepsilon} \int_{\R^2} \eta^2, \ \int_{\R^2} (\partial_1 N)^2 = \frac{18 \sqrt{2}}{\varepsilon^3} \int_{\R^2} (\partial_1 \eta)^2, \ \int_{\R^2} (\partial_2 N)^2 = \frac{36 \sqrt{2}}{\varepsilon^5} \int_{\R^2} (\partial_2 \eta)^2,$$
whereas
$$\int_{\R^2}(\partial_1 \Theta)^2 = \frac{36 \sqrt{2}}{\varepsilon} \int_{\R^2} (\partial_1 \varphi)^2, \ {\rm and} \ \int_{\R^2} (\partial_2 \Theta)^2 = \frac{72 \sqrt{2}}{\varepsilon^3} \int_{\R^2} (\partial_2 \varphi)^2.$$
A rather tedious computation along the same lines allows to derive the following expansions.

\begin{lemma}
\label{energyslow}
Let $v$ be a smooth map on $\R^2$ satisfying \eqref{loindezero}, and let $N$ and $\Theta$ be the corresponding functions defined by \eqref{slow-var}. The momentum $p(v)$ can be expressed in terms of the new functions as
\begin{equation}
\label{p-slow}
p(v) = \frac{\varepsilon}{72} \int_{\R^2} N \partial_1 \Theta,
\end{equation}
while the energy $E(v)$ has the expansion
\begin{equation}
\label{E-slow}
E(v) = \sqrt{2} \frac{\varepsilon}{144} \Big( E_0(N, \Theta) + \varepsilon^2 E_2(N, \Theta) + \varepsilon^4 E_ 4 (N, \Theta) \Big),
\end{equation}
where the functions $E_0$, $E_2$ and $E_4$ are given by
\begin{equation}
\label{defE0}
E_0(N, \Theta) = \int_{\R^2} \Big( N^2 + (\partial_1 \Theta)^2 \Big),
\end{equation}
\begin{equation}
\label{defE2}
E_2(N, \Theta) = \int_{\R^2} \Big( \frac{1}{2} (\partial_1 N)^2 + \frac{1}{2} (\partial_2 \Theta)^2 - \frac{1}{6} N (\partial_1 \Theta)^2 \Big),
\end{equation}
and
\begin{equation}
\label{defE4}
E_4(N, \Theta) = \int_{\R^2} \bigg( \frac{(\partial_2 N)^2}{4 - \frac{2 \varepsilon^2}{3} N} + \frac{N (\partial_1 N)^2}{12 - 2\varepsilon^2 N} - \frac{1}{12} N (\partial_2 \Theta)^2 \bigg).
\end{equation}
\end{lemma}

\begin{remark}
Recall that the map $u_\p$ found in Theorem \ref{dim2} minimizes the Ginzburg-Landau energy keeping the momentum $p$ fixed, equal to $\p$. If one takes instead only the first term of the energy in expansion \eqref{E-slow}, i.e. if one minimizes $E_0$ keeping the momentum $p$ equal to $\p$, then $\tilde{u}_\p$ will be a minimizer for the new problem if and only if
$$\tilde{N}_\p = \partial_1 \tilde{\Theta}_\p, \ {\rm and} \ \int_{\R^2} \tilde{N}_\p^2 = \frac{72 \p}{\varepsilon}.$$
Notice in particular that $\tilde{\Theta}_\p = \partial_1^{-1} \tilde{N}_\p$, so that $\partial_1^{-1} (\partial_2 \tilde{N}_\p) = \partial_2 \tilde{\Theta}_\p$. If we insert these relations into the definition of $E_2(\tilde{N}_\p, \tilde{\Theta}_\p)$, one obtains
$$E_2(\tilde{N}_\p, \tilde{\Theta}_\p) = \int_{\R^2} \bigg( \frac{1}{2} (\partial_1 \tilde{N}_\p)^2 + \frac{1}{2} (\partial_1^{-1} (\partial_2 \tilde{N}_\p))^2 - \frac{1}{6} \tilde{N}_\p^3 \bigg) = E_{KP}(\tilde{N}_\p).$$
This identity gives a first heuristic relation between the \eqref{GP} functional and the \eqref{KP} functional, as well as between the solutions $u_\p$ and the ground states for \eqref{KP}.
\end{remark}

Specifying the above change of variables to the case $v = u_\p$ and $\varepsilon = \varepsilon_\p$, setting $N_\p = N_{u_\p, \varepsilon_\p}$ and $\Theta_\p = \Theta_{u_\p, \varepsilon_\p}$, we obtain bounds for the integral quantities appearing in Lemma \ref{energyslow}. In view of \eqref{GLE} and \eqref{euhh}, we have
$$\int_{\R^2} (N_\p)^2 = \frac{18 \sqrt{2}}{\varepsilon_\p} \int_{\R^2} \eta_\p^2 \leq \frac{72 \sqrt{2} E(u_\p)}{\varepsilon_\p} \leq K,$$
where $K$ is some universal constant, whereas by \eqref{loindezero} and \eqref{relevenergy0},
$$\int_{\R^2} (\partial_1 \Theta_\p)^2 = \frac{36 \sqrt{2}}{\varepsilon_\p} \int_{\R^2} (\partial_1 \varphi_\p)^2 \leq \frac{144 \sqrt{2}}{\varepsilon_\p} \int_{\R^2} \varrho_\p^2 (\partial_1 \varphi_\p)^2 \leq \frac{288 \sqrt{2} E(u_\p)}{\varepsilon_\p},$$
so that
\begin{equation}
\label{primor}
\int_{\R^2} \bigg( (N_\p)^2 + (\partial_1 \Theta_\p)^2 \bigg) \leq K.
\end{equation}
Similarly, it follows from \eqref{deriveta} and \eqref{deriveta2} that
\begin{equation}
\label{dial}
\int_{\R^2} \bigg( (\partial_1 N_\p)^2 + (\partial_2 \Theta_\p)^2 \bigg) + \bigg| \int_{\R^2} N_\p (\partial_1 \Theta_\p)^2 \bigg| \leq K.
\end{equation}
For various other quantities, we only have at this stage rather crude estimates. For instance, concerning the uniform norm of $N_\p$, the bound provided by \eqref{elinfiniestimates} yields
\begin{equation}
\label{infinibound}
\| N_\p \|_{L^{\infty}(\R^2)} \leq \frac{K}{\varepsilon_\p^2}.
\end{equation}
We also only have for the transverse derivatives
\begin{equation}
\label{badestimate1}
\int_{\R^2} (\partial_2 N_\p)^2 + \int_{\R^2} \Big| N_\p (\partial_2 \Theta_\p)^2 \Big| \leq \frac{K}{\varepsilon_\p^2}.
\end{equation}
It follows from \eqref{dial} that
$$\big| E_2(N_\p, \Theta_\p) \big| \leq K,$$
whereas for $E_4$, we only obtain combining estimates \eqref{dial}, \eqref{infinibound} and \eqref{badestimate1},
$$\big| E_4(N_\p, \Theta_\p) \big| \leq \frac{K}{\varepsilon_\p^2}.$$
Hence, going back to the expansion of the energy, we deduce
\begin{equation}
\label{p3estim}
\Big| E(u_\p) - \sqrt{2} \frac{\varepsilon_\p}{144} E_0(N_\p, \Theta_\p) \Big| \leq K \varepsilon_\p^3.
\end{equation}
This leads to

\begin{lemma}
\label{difference}
There exists some positive constant $K$, not depending on $\p$, such that
\begin{equation}
\label{L2fort}
\int_{\R^2} \Big( N_\p - \partial_1 \Theta_\p \Big) ^2 \leq K \varepsilon_\p^2,
\end{equation}
for any $\p$ sufficiently small.
\end{lemma}

\begin{proof}
Using \eqref{p-slow}, \eqref{defE0} and \eqref{p3estim}, we are led to
$$\int_{\R^2} \Big( N_\p - \partial_1 \Theta_\p \Big) ^2 = E_0(N_\p, \Theta_\p) - 2 \int_{\R^2} N_\p \partial_1 \Theta_\p \leq \frac{144 E(u_\p)}{\sqrt{2} \varepsilon_\p} - \frac{144 \p}{\varepsilon_\p} + K \varepsilon_\p^2.$$
Since $E(u_\p) \leq \sqrt{2} \p$, the conclusion follows.
\end{proof}

Estimate \eqref{primor} provides a first step to compactness. In particular, there exists some map $N_0 \in L^2(\R^2)$ such that, up to a subsequence,
$$N_\p \rightharpoonup N_0 \ {\rm in} \ L^2(\R^2), \ {\rm as} \ \p \to 0.$$
As a consequence of Lemma \ref{difference}, we also have
$$\partial_1 \Theta_\p \rightharpoonup N_0 \ {\rm in} \ L^2(\R^2), \ {\rm as} \ \p \to 0.$$
To improve this convergence and characterize the limit function $N_0$, we turn to the equations for $N_\p$ and $\Theta_\p$.

\subsection{Expansion of the equations}

We now consider a finite energy solution $v$ to \eqref{TWc} satisfying \eqref{loindezero}, so that $v$ may be written as $v = \varrho \exp i \varphi$, and the functions $\varrho$ and $\varphi$ satisfy the system of equations \eqref{PolTWc1}-\eqref{PolTWc2}. At first order, each of the equations \eqref{PolTWc1} and \eqref{PolTWc2} express the fact that
$$N \sim \partial_1 \Theta, \ {\rm as} \ \varepsilon \to 0.$$
Indeed, we first have

\begin{lemma}
\label{slowequation1}
Assume $\varrho$ and $\varphi$ satisfy \eqref{PolTWc2}, and let $N$ and $\Theta$ be the corresponding functions defined by \eqref{slow-var}. Then, $N$ and $\Theta$ satisfy
\begin{equation}
\label{slow-1}
N - \partial_1 \Theta = \varepsilon^2 \Big( \boL_{\varepsilon, 1}(N, \Theta) + \boR_{\varepsilon, 1}(N, \Theta) \Big), 
\end{equation}
where the remainder terms $\boL_{\varepsilon, 1}(N, \Theta)$ and $\boR_{\varepsilon, 1}(N, \Theta)$ are given by
$$\boL_{\varepsilon, 1}(N, \Theta) = \frac{1}{\varepsilon^2} \Big( \sqrt{1 - \frac{\varepsilon^2}{2}} - 1 \Big) \partial_1 \Theta + \frac{1}{2} \partial_1^2 N + \frac{\varepsilon^2}{4} \partial_2^2 N,$$
and
\begin{align*}
\begin{split}
\boR_{\varepsilon,1}(N, \Theta) = & \frac{1}{12} \bigg( 2 N^2 - 2 \sqrt{1 - \frac{\varepsilon^2}{2}} N \partial_1 \Theta + (\partial_1 \Theta)^2 \bigg)\\
+ & \frac{\varepsilon^2}{72} \bigg( 3 \frac{(\partial_1 N)^2}{1 - \frac{\varepsilon^2}{6} N} - N (\partial_1 \Theta)^2 + 3 (\partial_2 \Theta)^2 \bigg)\\
+ & \frac{\varepsilon^4}{144} \bigg( 3 \frac{(\partial_2 N)^2}{1 - \frac{\varepsilon^2}{6} N} - N (\partial_2\Theta)^2 \bigg).
\end{split}
\end{align*}
\end{lemma}

We similarly have

\begin{lemma}
\label{slowequation2}
Assume $\varrho$ and $\varphi$ satisfy \eqref{PolTWc1}, and let $N$ and $\Theta$ be the corresponding functions defined by \eqref{slow-var}. Then, $N$ and $\Theta$ satisfy
\begin{equation}
\label{slow-2}
\partial_1 N - \partial_1^2 \Theta = \varepsilon^2 \Big( \boL_{\varepsilon, 2}(N, \Theta) + \boR_{\varepsilon, 2}(N, \Theta) \Big),
\end{equation}
where the remainder terms $\boL_{\varepsilon, 2}(N, \Theta)$ and $\boR_{\varepsilon, 2}(N, \Theta) $ are given by
$$\boL_{\varepsilon, 2}(N, \Theta) = \frac{1}{\varepsilon^2} \Big( 1 - \sqrt{1 - \frac{\varepsilon^2}{2}} \Big) \partial_1 N + \frac{1}{2} \partial_2^2 \Theta,$$
and
$$\boR_{\varepsilon, 2}(N, \Theta) = - \frac{1}{6} \partial_1 \big[ N \partial_1 \Theta \big] - \frac{\varepsilon^2}{12} \partial_2 \big[ N \partial_2 \Theta \big].$$
\end{lemma}

As mentioned above, equations \eqref{slow-1} and \eqref{slow-2} twice express the fact that the functions $N$ and $\partial_1 \Theta$ are equal at the limit $\varepsilon \to 0$. In order to identify their common limit, we expand some combination of \eqref{slow-1} and \eqref{slow-2} to deduce

\begin{prop}
\label{prop-slow-N}
Let $v$ be a finite energy solution to \eqref{TWc} on $\R^2$ satisfying \eqref{loindezero}, and let $N$ and $\Theta$ be the corresponding functions defined by \eqref{slow-var}. Then, $N$ and $\Theta$ satisfy
\begin{equation}
\label{N-eq}
\boL (N) = - \partial_1^2 \Big[ \frac{1}{3} N^2 + \frac{1}{6} (\partial_1 \Theta)^2 \Big] + \varepsilon^2 \Big( \boL_\varepsilon(N) + \boR_\varepsilon(N, \Theta) \Big), 
\end{equation}
where $\boL$ is the linear operator given by
$$\boL (N) = \partial_1^4 N - \Delta N,$$
and the remainder terms $\boL_\varepsilon(N)$ and $\boR_\varepsilon(N, \Theta)$ are given by
$$\boL_\varepsilon(N) = - \partial_1^2 \partial_2^2 N - \frac{\varepsilon^2}{4}\partial_2^4 N,$$
and
\begin{align*}
\boR_\varepsilon(N, \Theta) = & \frac{1}{72} \bigg( 2 \partial_1^2 \big[ N (\partial_1 \Theta)^2 \big] - 6 \partial_1^2 \Big[\frac{(\partial_1 N)^2}{1 - \frac{\varepsilon^2}{6} N} \Big] - 24\partial_2^2 (N^2)- 6 \sqrt{1 - \frac{\varepsilon^2}{2}} \partial_1 \partial_2 \big[ N\partial_2 \Theta \big]\\
+ & 12 \sqrt{1 - \frac{\varepsilon^2}{2}}\partial_2^2 \big[ N \partial_1 \Theta \big] - 3 \partial_1^2 \big[ (\partial_2 \Theta)^2 \big] - 6\partial_2^2 \big[ (\partial_1 \Theta)^2 \big] \bigg)\\
+ & \frac{\varepsilon^2}{144} \bigg( - 3 \partial_1^2 \Big[ \frac{(\partial_2 N)^2}{1 - \frac{\varepsilon^2}{6} N} \Big] + \partial_1^2 \big[ N (\partial_2 \Theta)^2 \big] - 6 \partial_2^2 \Big[\frac{(\partial_1 N)^2}{1 - \frac{\varepsilon^2}{6} N} \Big] + 2\partial_2^2 \big[ N (\partial_1 \Theta)^2 \big]\\
- & 3 \partial_2^2 \big[ (\partial_2 \Theta)^2 \big] \bigg) + \frac{\varepsilon^4}{288} \bigg( - 3 \partial_2^2 \Big[ \frac{(\partial_2 N)^2}{1 - \frac{\varepsilon^2}{6} N} \Big] + \partial_2^2 \big[ N (\partial_2 \Theta)^2\big] \bigg).
\end{align*}
\end{prop}

\begin{proof}
Equation \eqref{N-eq} is derived applying the differential operator $- \partial_1^2 - \frac{\varepsilon^2}{2} \partial_2^2$ to \eqref{slow-1}, the operator $\sqrt{1 - \frac{\varepsilon^2}{2}} \partial_1$ to \eqref{slow-2}, and adding the corresponding relations.
\end{proof}

Notice that we have at this stage,
$$\partial_1^4 N - \Delta N + \frac{1}{2}\partial_1^2 N^2 = \frac{1}{6} \partial_1^2 (N^2 - (\partial_1 \Theta)^2) + \varepsilon^2 \Big( \boL_\varepsilon(N, \Theta) + \boR_\varepsilon(N, \Theta) \Big),$$
where we recognize equation \eqref{eqSW2} for $N$ in the left-hand side. Specifying this relation to the solutions $N_\p$ and $\Theta_\p$, it remains to prove that the weak limit $N_0$ of the sequence $(N_\p)_{\p > 0}$ is a solution to \eqref{SW}, and to show some strong convergence. This requires to establish that the nonlinear remainder term $\boR_\varepsilon$ is small in some suitable sense. Indeed, the first term on the right-hand side will tend to $0$ in view of Lemma \ref{difference}, whereas the linear term $\boL_\varepsilon(N)$ presents no difficulty.

The remainder term $\boR_\varepsilon$ is a sum of several second order derivatives. We order them according to the type of second order derivatives, writing
$$\boR_\varepsilon(N, \Theta) = \sum_{i + j = 2} \partial_1^i \partial_2^j \boR_\varepsilon^{i, j},$$
where
\begin{equation}
\label{r20}
\boR_\varepsilon^{2, 0} = \frac{1}{36} N (\partial_1 \Theta)^2 - \frac{(\partial_1 N)^2}{12 (1 - \frac{\varepsilon^2}{6} N)} - \frac{1}{24} (\partial_2 \Theta)^2 - \varepsilon^2 \frac{(\partial_2 N)^2}{48 (1 - \frac{\varepsilon^2}{6} N)} + \frac{\varepsilon^2}{144} N (\partial_2 \Theta)^2,
\end{equation}
\begin{equation}
\label{r02}
\begin{split}
\boR_\varepsilon^{0, 2} & = - \frac{N^2}{3} + \frac{1}{6}\sqrt{1 - \frac{\varepsilon^2}{2}} N \partial_1 \Theta - \frac{(\partial_1 \Theta)^2}{12} - \varepsilon^2 \frac{(\partial_1 N)^2}{24 (1 - \frac{\varepsilon^2}{6} N)} + \frac{\varepsilon^2}{72} N (\partial_1 \Theta)^2 - \frac{\varepsilon^2}{48} (\partial_2 \Theta)^2\\
& - \varepsilon^4 \frac{(\partial_2 N)^2}{96 (1 - \frac{\varepsilon^2}{6}N)} + \frac{\varepsilon^4}{288} N (\partial_2 \Theta)^2,\\
\end{split}
\end{equation}
and
\begin{equation}
\label{r11}
\boR_\varepsilon^{1, 1} = - \frac{1}{12} \sqrt{1 - \frac{\varepsilon^2}{2}} N \partial_2 \Theta.
\end{equation}

In several places, it will be convenient to write
$$\boR^{i, j}_\varepsilon = R^{i, j}_\varepsilon + \varepsilon^2 \nu^{i, j}_\varepsilon,$$
where $\nu^{1, 1}_\varepsilon = 0$,
$$\nu^{2, 0}_\varepsilon = - \frac{(\partial_2 N)^2}{48 (1 - \frac{\varepsilon^2}{6} N)}+ \frac{1}{144} N (\partial_2 \Theta)^2,$$
and
$$\nu^{0, 2}_\varepsilon = - \frac{(\partial_1 N)^2}{24 (1 - \frac{\varepsilon^2}{6} N)} + \frac{1}{72} N (\partial_1 \Theta)^2 - \frac{1}{48}(\partial_2 \Theta)^2 - \varepsilon^2 \frac{(\partial_2 N)^2}{96 (1 - \frac{\varepsilon^2}{6} N)} + \frac{\varepsilon^2}{288} N (\partial_2 \Theta)^2.$$
Notice in particular that
\begin{equation}
\label{anquetil}
|R^{2, 0}_\varepsilon| \leq K \Big( (\partial_1 N)^2 + (\partial_2 \Theta)^2 + |N (\partial_1 \Theta)^2| \Big),
\end{equation}
whereas
\begin{equation}
\label{poulidor}
|R^{1, 1}_\varepsilon| \leq K \big| N \big| \big| \partial_2 \Theta \big|, \ {\rm and} \ |R^{0, 2}_\varepsilon| \leq K \Big( N^2 + (\partial_1 \Theta)^2 \Big).
\end{equation}
Similarly, we also have
\begin{equation}
\label{merckx}
\begin{split}
|\nu^{2, 0}_\varepsilon| & \leq K \Big( (\partial_2 N)^2 + |N (\partial_2 \Theta)^2| \Big),\\
|\nu^{0, 2}_\varepsilon| & \leq K \Big( (\partial_1 N)^2 + |N (\partial_1 \Theta)^2| + (\partial_2 \Theta)^2 + \varepsilon^2 \big( (\partial_2 N)^2 + |N (\partial_2 \Theta)^2| \big) \Big).
\end{split}
\end{equation}
Specifying the previous quantities for $N_\p$ and $\Theta_\p$, we obtain some initial bounds on the nonlinear remainder terms, which will prove essential to compute the estimates of Proposition \ref{Sobolevbound}.

\begin{lemma}
\label{cosa}
There exists some positive constant $K$, not depending on $\p$, such that
\begin{equation}
\label{cosa1}
\int_{\R^2} \Big( |\boR_{\varepsilon_\p}^{1, 1}| + |R_{\varepsilon_\p}^{0, 2}| \Big) \leq K,
\end{equation}
and
\begin{equation}
\label{cosa2}
\int_{\R^2} \Big( |R_{\varepsilon_\p}^{2, 0}| + |\nu^{2, 0}_{\varepsilon_\p}| + |\nu^{0, 2}_{\varepsilon_\p}| \Big) \leq \frac{K}{\varepsilon_\p^2},
\end{equation}
for any $\p$ sufficiently small.
\end{lemma}

\begin{proof}
Bounds \eqref{cosa1} and \eqref{cosa2} are consequences of bounds
\eqref{primor}, \eqref{dial} and \eqref{badestimate1}, and inequalities \eqref{anquetil}, \eqref{poulidor} and \eqref{merckx}. Concerning the term $\int_{\R^2} N_\p (\partial_1 \Theta_\p)^2$ in \eqref{cosa2}, we have to invoke the crude bound \eqref{infinibound}, which yields
$$\int_{\R^2} \Big| N_\p (\partial_1 \Theta_\p)^2 \Big| \leq \frac{K}{\varepsilon_\p^2} \int_{\R^2} (\partial_1 \Theta_\p)^2 \leq \frac{K}{\varepsilon_\p^2}.$$
\end{proof}

\subsection{Estimates for the phase $\Theta_\p$}
\label{Tetee}

In the previous discussion, we did not consider the function $\Theta$. In particular, we did not compute any rescaled equation for this function. Applying the partial differential operator $\boL - \varepsilon^2 \boL_\varepsilon$ to \eqref{slow-2} and introducing equation \eqref{N-eq} in the resulting equation in order to eliminate the function $N$ in the linear part, we compute
\begin{equation}
\label{Theta-eq}
\boL (\partial_1^2 \Theta) = - \partial_1^3 \Big( \frac{1}{3} N^2 + \frac{1}{6} (\partial_1 \Theta)^2 \Big) + \varepsilon^2 \Big( \boL_{\varepsilon, 3}(\Theta) + \boR_{\varepsilon, 3}(N, \Theta) \Big),
\end{equation}
where the remainder terms $\boL_{\varepsilon, 3}(\Theta)$ and $\boR_{\varepsilon, 3}(N, \Theta)$ are given by
$$\boL_{\varepsilon, 3}(\Theta) = \boL_\varepsilon(\partial_1^2 \Theta) - \frac{1}{2} \boL(\partial_2^2 \Theta) + \frac{\varepsilon^2}{2} \boL_\varepsilon (\partial_2^2 \Theta),$$
and
\begin{align*}
\boR_{\varepsilon, 3}(N, \Theta) & = \frac{1}{\varepsilon^2} \bigg( 1 - \sqrt{1 - \frac{\varepsilon^2}{2}} \bigg) \partial_1^3 \Big( \frac{1}{3} N^2 + \frac{1}{6} (\partial_1 \Theta)^2 \Big) + \sqrt{1 - \frac{\varepsilon^2}{2}} \partial_1 \boR_\varepsilon(N, \Theta)\\
& - \boL \big( \boR_{\varepsilon, 2}(N, \Theta) \big) + \varepsilon^2 \boL_\varepsilon \big( \boR_{\varepsilon, 2}(N, \Theta) \big).
\end{align*}
At least formally, this may be written as
$$\partial_1^4 (\partial_1 \Theta) - \Delta(\partial_1 \Theta) + \frac{1}{2} \partial_1^2 (\partial_1 \Theta)^2 = \frac{1}{3} \partial_1^2 \Big( (\partial_1 \Theta)^2 - N^2 \Big) + \varepsilon^2 \partial_1^{- 1} \Big( \boL_{\varepsilon, 3}(\Theta) + \boR_{\varepsilon, 3}(N, \Theta) \Big).$$
We recognize once more equation \eqref{eqSW2} for $\partial_1 \Theta$ in the left-hand side. However, the analysis of equation \eqref{Theta-eq} is substantially more difficult than the study of \eqref{N-eq}, due to the intricacy of the remainder terms and the necessity to apply the operator $\partial_1^{- 1}$ to \eqref{Theta-eq} to recover \eqref{eqSW2}. Hence, our argument to deal with the phase $\Theta$ does not rely on \eqref{Theta-eq}. Instead, we invoke the estimates of Lemma \ref{Derivphi}, whose rescaled versions give bounds on $\Theta$ in function of those on $N$.

\begin{lemma}
\label{Scaledtheta}
Let $1 < q < + \infty$. There exists some positive constant $K(q)$, not depending on $\p$, such that
\begin{equation}
\label{thetap}
\| \partial_1 \Theta_\p \|_{L^q(\R^2)} + \varepsilon_\p \| \partial_2 \Theta_\p \|_{L^q(\R^2)} \leq K(q) \| N_\p \|_{L^q(\R^2)},
\end{equation}
for any $\p$ sufficiently small. Similarly, given any $\alpha \in \N^2$, and denoting
$$\Xi_\p(q, \alpha) \equiv \| \partial^\alpha \partial_1 \Theta_\p \|_{L^q(\R^2)} + \varepsilon_\p \| \partial^\alpha \partial_2 \Theta_\p \|_{L^q(\R^2)},$$
there exists some positive constant $K(q, \alpha)$, not depending on $\p$, such that
\begin{equation}
\label{thepat}
\Xi_\p(q, \alpha) \leq K(q, \alpha) \bigg( \| \partial^\alpha N_\p \|_{L^q(\R^2)} + \varepsilon_\p^2 \sum_{0 \leq \beta < \alpha} \| \partial^\beta N_\p \|_{L^\infty(\R^2)} \Xi_\p(q, \alpha - \beta) \bigg),
\end{equation}
for any $\p$ sufficiently small.
\end{lemma}

\begin{proof}
Inequalities \eqref{thetap} and \eqref{thepat} are rescaled versions of \eqref{pepe} and \eqref{meme}. In view of scalings \eqref{slow-var}, given any $1 < q \leq + \infty$, the $L^q$-norm of the function $\partial^\alpha N$ is related to the $L^q$-norm of $\partial^\alpha \eta$ by
\begin{equation}
\label{mauresmo}
\| \partial^\alpha N \|_{L^q(\R^2)} = \frac{K(q, \alpha)}{\varepsilon^{2 + \alpha_1 + 2 \alpha_2 - \frac{3}{q}}} \| \partial^\alpha \eta \|_{L^q(\R^2)},
\end{equation}
where $K(q, \alpha)$ denotes some positive constant, not depending on $\varepsilon$. Similarly, we compute for the functions $\partial^\alpha \partial_1 \Theta$ and $\partial^\alpha \partial_2 \Theta$,
\begin{equation}
\label{aencoreloose}
\| \partial^\alpha \partial_1 \Theta \|_{L^q(\R^2)} = \frac{K(q, \alpha)}{\varepsilon^{2 + \alpha_1 + 2 \alpha_2 - \frac{3}{q}}} \| \partial^\alpha \partial_1 \varphi \|_{L^q(\R^2)}, \ {\rm and} \ \| \partial^\alpha \partial_2 \Theta \|_{L^q(\R^2)} = \frac{K(q, \alpha)}{\varepsilon^{3 + \alpha_1 + 2 \alpha_2 - \frac{3}{q}}} \| \partial^\alpha \varphi \|_{L^q(\R^2)}.
\end{equation}
Inequalities \eqref{thetap} and \eqref{thepat} then follow from rescaling \eqref{pepe} and \eqref{meme}, specifying identities \eqref{mauresmo} and \eqref{aencoreloose} for the functions $N_\p$ and $\Theta_\p$.
\end{proof}

In view of Lemma \ref{Scaledtheta}, we will not invoke equation \eqref{Theta-eq} to bound the function $\Theta_\p$. Instead, we will take advantage of the regularizing properties of equation \eqref{N-eq}, and rely on the initial estimates of Lemma \ref{cosa}, to bound the $L^q$-norm of $N_\p$ (and actually, its first order derivatives) independently on $\p$. We will then deduce from \eqref{thetap} and \eqref{thepat}, $L^q$-estimates of some low order derivatives of $\Theta_\p$. This in turn will provide new bounds on the nonlinear terms $\boR_{\varepsilon_\p}^{i,j}$, and on their first order derivatives, improving the estimates of Lemma \ref{cosa}. Using in particular, the inductive nature of \eqref{thepat}, we will iterate the argument to obtain $L^q$-bounds on any order derivatives of $N_\p$ and $\Theta_\p$, and complete the proof of Proposition \ref{Sobolevbound} (see Section \ref{Unity} below). Notice that this strategy will first require to analyse the regularizing nature of \eqref{N-eq} which becomes more transparent taking its Fourier transform.

\subsection{Kernels of the rescaled equations}
\label{Abricot}

We derive a new formulation of \eqref{N-eq} which brings out its regularizing properties. Taking the Fourier transform of the previous rescaled equations, we deduce

\begin{cor}
\label{cor-fou}
Let $v$ be a finite energy solution to \eqref{TWc} on $\R^2$ satisfying \eqref{loindezero}, and let $N$ and $\Theta$ be the corresponding functions defined by \eqref{slow-var}. Then, $\widehat{N}$ and $\widehat{\Theta}$ satisfy
\begin{equation}
\label{slow-fou-1}
\Big( 1 + \frac{\varepsilon^2}{2} \xi_1^2 + \frac{\varepsilon^4}{4} \xi_2^2 \Big) \widehat{N}(\xi) - i \sqrt{1 - \frac{\varepsilon^2}{2}} \xi_1 \widehat{\Theta}(\xi) = \varepsilon^2 \widehat{\boR_{\varepsilon, 1}}(\xi),
\end{equation}
\begin{equation}
\label{slow-fou-2}
\Big( \xi_1^2 + \frac{\varepsilon^2}{2} \xi_2^2 \Big) \widehat{\Theta}(\xi) + i \sqrt{1 - \frac{\varepsilon^2}{2}} \xi_1 \widehat{N}(\xi) = \varepsilon^2 \widehat{\boR_{\varepsilon, 2}}(\xi),
\end{equation}
and
\begin{equation}
\label{N-fou-eq}
\Big( \xi_1^4 + |\xi|^2 + \varepsilon^2 \xi_1^2 \xi_2^2 + \frac{\varepsilon^4}{4} \xi_2^4 \Big) \widehat{N}(\xi) = \xi_1^2 \Big( \frac{1}{3} \widehat{N^2}(\xi) + \frac{1}{6} \widehat{(\partial_1 \Theta)^2}(\xi) \Big) + \varepsilon^2 \widehat{\boR_\varepsilon}(\xi).
\end{equation}
\end{cor}

\begin{proof}
Equations \eqref{slow-fou-1}, \eqref{slow-fou-2} and \eqref{N-fou-eq} follow from taking the Fourier transform of equations \eqref{slow-1}, \eqref{slow-2} and \eqref{N-eq}.
\end{proof}

At this stage, it is presumably worthwhile to compare equations \eqref{N-fou-eq} and \eqref{fouSW}. This leads us to consider the perturbed kernel $K_\varepsilon$, whose Fourier transform is given by
$$\widehat{K_\varepsilon}(\xi) = \frac{\xi_1^2}{|\xi|^2 + \xi_1^4 + \varepsilon^2 \xi_1^2 \xi_2^2 + \frac{\varepsilon^4}{4} \xi_2^4}.$$
The kernel $K_\varepsilon$ is a regularization of the kernel $K_0$, since it belongs to $H^\frac{1}{4}(\R^2)$ (see Lemma \ref{noyaudur} below), and tends to $K_0$ in $L^2(\R^2)$, as $\varepsilon \to 0$, by the dominated convergence theorem. We will extensively use this additional regularizing property of $K_\varepsilon$ to compute estimates of the function $N$.

More generally, since
$$\widehat{\boR_\varepsilon}(\xi) = - \sum_{i + j = 2} \xi_i \xi_j \widehat{\boR_\varepsilon}^{i, j}(\xi),$$
we also introduce the kernels $K_\varepsilon^{i, j}$ defined by
\begin{equation}
\label{defkernel}
\widehat{K_\varepsilon^{i, j}}(\xi) = \frac{\xi_1^i \xi_2^j}{|\xi|^2 + \xi_1^4 + \varepsilon^2 \xi_1^2 \xi_2^2 + \frac{\varepsilon^4}{4} \xi_2^4},
\end{equation}
for any $0 \leq i, j \leq 4$ such that $2 \leq i + j \leq 4$ (so that, in particular, $K_\varepsilon = K_\varepsilon^{2, 0}$). We then recast equation \eqref{N-eq} as a convolution equation
\begin{equation}
\label{decomposition}
N_\p = K^{2, 0}_{\varepsilon_\p} \star f_\p - \sum_{i + j = 2} \varepsilon_\p^2 K^{i,j}_{\varepsilon_\p} \star \boR^{i, j}_{\varepsilon_\p},
\end{equation}
where
\begin{equation}
\label{effet}
f_\p = \frac{1}{3} N_\p^2 + \frac{1}{6} (\partial_1 \Theta_\p)^2.
\end{equation}
In view of the multiplier properties of the kernels $K_{\varepsilon_\p}^{i, j}$ (see Lemma \ref{Multiprop} below), equation \eqref{decomposition} provides a control on the $L^q$-norm of $N_\p$ in function of the $L^q$-norms of the nonlinear terms $f_\p$ and $\boR_{\varepsilon_\p}^{i,j}$. This control is the starting point of the proof of Proposition \ref{Sobolevbound}, which follows combining the superlinear nature of the nonlinear terms $f_\p$ and $\boR_{\varepsilon_\p}^{i,j}$ with the estimates of $\Theta_\p$ provided by Lemma \ref{Scaledtheta} (see Section \ref{Unity} below).

\section{Properties of the kernels $K_\varepsilon^{i, j}$}
\label{Linear}

We now turn to the analysis of the kernels $K_\varepsilon^{i, j}$. In particular, we provide a number of estimates, which are required by the proof of Proposition \ref{Sobolevbound}.

\subsection{$H^\alpha$-estimates of the kernels}

For given $0 \leq \alpha < 1$, we establish $H^\alpha$-estimates for the kernels $K_\varepsilon^{i, j}$. We first consider their $\dot{H}^\alpha$-semi-norms defined in the Fourier space by
$$\| K^{i, j}_\varepsilon \|_{\dot{H}^\alpha(\R^2)}^2 = \int_{\R^2} |\xi|^{2 \alpha} |\widehat{K_\varepsilon^{i, j}}(\xi)|^2 d\xi.$$

\begin{lemma}
\label{noyaudur}
Let $0 < \varepsilon \leq 1$ and $0 \leq \alpha < 1$. Then,
\begin{equation}
\label{keke1}
\| K^{2, 0}_\varepsilon \|_{\dot{H}^\alpha(\R^2)} \leq K(\alpha) \big( 1 + \varepsilon^{\frac{1}{2} - 2 \alpha} \big), \ \| K^{1, 1}_\varepsilon \|_{\dot{H}^\alpha(\R^2)} \leq K(\alpha) \big( 1 + \varepsilon^{- \frac{1}{2} - 2 \alpha} \big),
\end{equation}
and
\begin{equation}
\label{keke2}
\| K^{0, 2}_\varepsilon \|_{\dot{H}^\alpha(\R^2)} \leq K(\alpha) \big( 1 + \varepsilon^{- \frac{3}{2} - 2 \alpha} \big).
\end{equation}
\end{lemma}

\begin{proof}
The proof is an explicit computation. In view of the definition of the semi-norms, we compute using polar coordinates, and noticing that $i + j = 2$,
\begin{align*}
\| K^{i, j}_\varepsilon \|_{\dot{H}^\alpha(\R^2)}^2 & = \int_{\R^2} \frac{|\xi|^{2 \alpha} \xi_1^{2 i} \xi_2^{2 j}}{\big( |\xi|^2 + \xi_1^4 + \varepsilon^2 \xi_1^2 \xi_2^2 + \frac{\varepsilon^4}{4} \xi_2^4 \big)^2} d\xi\\
& = \int_0^{+ \infty} \int_0^{2 \pi} r^{2 \alpha + 1} \frac{\cos (\theta)^{2 i} \sin(\theta)^{2 j}}{\big( 1 + r^2 \cos(\theta)^4 + \varepsilon^2 r^2 \cos(\theta)^2 \sin(\theta)^2 + \frac{\varepsilon^4}{4} r^2 \sin(\theta)^4 \big)^2} dr d\theta\\
& = 4 \int_0^{+ \infty} \int_0^{+ \infty} r^{2 \alpha + 1} \frac{u^{2 j} (1 + u^2)^{3 - i - j}}{\big( (1 + u^2)^2 + r^2 + \varepsilon^2 r^2 u^2 + \frac{\varepsilon^4}{4} r^2 u^4 \big)^2} dr du,
\end{align*}
where we have set $u = \tan(\theta)$ in the last integral. The previous computation leads us to introduce the quantity
$$J_{\beta, \varepsilon}(r) = \int_0^{+ \infty} \frac{u^{2 \beta}}
{\big( (1 + u^2)^2 + r^2 + \varepsilon^2 r^2 u^2 + \frac{\varepsilon^4}{4} r^2 u^4 \big)^2} du,$$
so that
\begin{equation}
\label{noyaute}
\| K^{i, j}_\varepsilon \|_{\dot{H}^\alpha(\R^2)}^2 \leq K \int_0^{+ \infty} r^{2 \alpha + 1} \big( J_{\beta_1, \varepsilon}(r) + J_{\beta_2, \varepsilon}(r) \big) dr,
\end{equation}
where $\beta_1 = j$ and $\beta_2 = 3 - i$. We now claim that
\begin{equation}
\label{reclame}
\int_0^{+ \infty} r^{2 \alpha + 1} J_{\beta, \varepsilon}(r) dr \leq K(\alpha, \beta) \bigg( 1 + \frac{1}{\varepsilon^{4 \alpha + 2 \beta - 3}} \bigg),
\end{equation}
for any $0 \leq \beta < \frac{7}{2}$ and any $0 \leq \alpha < 1$. We postpone the proof of Claim \eqref{reclame}, and first complete the proof of Lemma \ref{noyaudur}. Combining identity \eqref{noyaute} with \eqref{reclame}, we obtain
$$\| K^{i, j}_\varepsilon \|_{\dot{H}^\alpha(\R^2)}^2 \leq K
\bigg( 1 + \frac{1}{\varepsilon^{4 \alpha + 2 j - 3}} + \frac{1}{\varepsilon^{4 \alpha + 3 - 2 i}} \bigg),$$
and the conclusion follows applying this inequality for the various choices of $i$ and $j$.
\end{proof}

\begin{proof}[Proof of Claim \eqref{reclame}]
In order to estimate the integral in the left-hand side of Claim \eqref{reclame}, we first compute some bounds for the function $J_{\beta, \varepsilon}$. When $0 \leq r \leq 1$, we have
\begin{equation}
\label{jiji1}
|J_{\beta, \varepsilon}(r)| \leq \int_0^{+ \infty} \frac{u^{2 \beta}}{1 + u^8} du \leq K(\beta),
\end{equation}
since $0 \leq \beta < \frac{7}{2}$. On the other hand, when $r > 1$, we compute
$$|J_{\beta, \varepsilon}(r)| \leq K \Bigg( \int_0^1 \frac{du}{1 + r^4} + \int_1^\frac{1}{\varepsilon} \frac{u ^{2 \beta}}{(u^4 + r^2)^2} du + \int_\frac{1}{\varepsilon}^{+ \infty} \frac{u ^{2 \beta - 8}}{(1 + r^2 \varepsilon^4)^2} du \Bigg),$$
so that, since $0 \leq \beta < \frac{7}{2}$,
\begin{equation}
\label{jiji2}
|J_{\beta, \varepsilon}(r)| \leq K(\beta) \Big( \frac{1}{r^4} + r^{\beta - \frac{7}{2}} + \varepsilon^{7 - 2 \beta} \Big),
\end{equation}
when $1 \leq r \leq \frac{1}{\varepsilon^2}$. Similarly, when $r \geq \frac{1}{\varepsilon^2}$,
\begin{equation}
\label{jiji3}
|J_{\beta, \varepsilon}(r)| \leq K(\beta) \Big( \frac{1}{r^4} + \frac{1}{\varepsilon^{2 \beta + 1} r^4} + \frac{\varepsilon^{7 - 2 \beta}}{(1 + r^2 \varepsilon^4)^2} \Big).
\end{equation}
Estimates \eqref{jiji1}, \eqref{jiji2} and \eqref{jiji3} finally provide Claim \eqref{reclame}, when $0 \leq \alpha < 1$.
\end{proof}

Since inequalities \eqref{keke1} and \eqref{keke2} are also valid for $\alpha = 0$, i.e. for the $L^2$-norm, we may remove the dots in inequalities \eqref{keke1} and \eqref{keke2}. Notice in particular that we have the bounds
\begin{equation}
\label{crepewahou}
\| K^{1, 1}_\varepsilon \|_{H^\alpha(\R^2)} + \varepsilon \| K^{1, 2}_\varepsilon \|_{H^\alpha(\R^2)} + \varepsilon^2 \| K^{2, 2}_\varepsilon \|_{H^\alpha(\R^2)} \leq K(\alpha),
\end{equation}
for any $0 \leq \alpha \leq \frac{1}{4}$.

\subsection{Multiplier properties of the kernels}

We now provide some multiplier properties of the kernels $K^{i, j}_\varepsilon$. Our analysis relies on a theorem by Lizorkin \cite{Lizorki1}
\footnote{Estimate \eqref{ibanez} in Theorem \ref{thelizo} is more precisely a consequence of Lemma 6 and of the proof of Theorem 8 in \cite{Lizorki1}.}
, which we first recall for sake of completeness.

\begin{theorem}[\cite{Lizorki1}]
\label{thelizo}
Let $\widehat{K}$ be a bounded function in $\boC^2(\R^2 \setminus \{ 0 \})$, and assume that
$$\xi_1^{k_1} \xi_2^{k_2} \partial_1^{k_1} \partial_2^{k_2} \widehat{K}(\xi) \in L^\infty(\R^2),$$
for any $0 \leq k_1, k_2 \leq 1$ such that $k_1 + k_2 \leq 2$. Then, $\widehat{K}$ is a multiplier from $L^q(\R^2)$ to $L^{q}(\R^2)$ for any $1 < q < + \infty$. More precisely, given any $1 < q < + \infty$, there exists a constant $K(q)$, depending only on $q$, such that
\begin{equation}
\label{ibanez}
\| K \star f \|_{L^q(\R^2)} \leq K(q) M(\widehat{K}) \| f \|_{L^q(\R^2)}, \ \forall f \in L^q(\R^2),
\end{equation}
where we denote
$$M(\widehat{K}) \equiv \sup \Big\{ |\xi_1|^{k_1} |\xi_2|^{k_2} \Big| \partial_1^{k_1} \partial_2^{k_2} \widehat{K}(\xi) \Big|, \xi \in \R^2, 0 \leq k_1 \leq 1, 0 \leq k_2 \leq 1, k_1 + k_2 \leq 2 \Big\}.$$
\end{theorem}

Applying Theorem \ref{thelizo} to the kernels $K^{i, j}_\varepsilon$, we obtain

\begin{lemma}
\label{Multiprop}
Let $1 < q < + \infty$. Given any integers $0 \leq i, j \leq 4$ such that $2 \leq i + j \leq 4$, we denote
$$\kappa_{i, j} = \max \{ i + 2 j - 4, 0 \},$$
Then, there exists some positive constant $K(q)$, not depending on $\varepsilon$, such that
\begin{equation}
\label{multiLq}
\| K^{i, j}_\varepsilon \star f \|_{L^q(\R^2)} \leq \frac{K(q)}{\varepsilon^{\kappa_{i ,j}}} \| f \|_{L^q(\R^2)},
\end{equation}
for any function $f \in L^q(\R^2)$ and any $\varepsilon > 0$.
\end{lemma}

\begin{proof}
Inequality \eqref{multiLq} is a consequence of \eqref{ibanez} once we have checked that the functions $\widehat{K^{i, j}_\varepsilon}$ satisfy the assumptions of Theorem \ref{thelizo}, and established the dependence with respect to $\varepsilon$ of the quantity $M(\widehat{K_\varepsilon^{i, j}})$.

First notice that the functions $\widehat{K^{i, j}_\varepsilon}$, which are bounded on $\R^2$, and belong to $\boC^2(\R^2 \setminus \{ 0 \})$, may be written as
$$\widehat{K_\varepsilon^{i, j}}(\xi) = \frac{\xi_1^i \xi_2^j}{Q(\xi)},$$
where $Q(\xi) \equiv |\xi|^2 + \xi_1^4 + \varepsilon^2 \xi_1^2 \xi_2^2 + \frac{\varepsilon^4}{4} \xi_2^4$. We therefore compute
\begin{equation}
\label{namur}
\xi_1 \partial_1 \widehat{K_\varepsilon^{i, j}}(\xi) = i \frac{\xi_1^i \xi_2^j}{Q(\xi)} - \frac{\xi_1^i \xi_2^j}{Q(\xi)} \frac{\xi_1 \partial_1 Q(\xi)}{Q(\xi)}, \ \xi_2 \partial_2 \widehat{K_\varepsilon^{i, j}}(\xi) = j \frac{\xi_1^i \xi_2^j}{Q(\xi)} - \frac{\xi_1^i \xi_2^j}{Q(\xi)} \frac{\xi_2 \partial_2 Q(\xi)}{Q(\xi)},
\end{equation}
and
\begin{equation}
\label{charleroi}
\xi_1 \xi_2 \partial_1 \partial_2 \widehat{K_\varepsilon^{i, j}}(\xi) = \frac{\xi_1^i \xi_2^j}{Q(\xi)} \bigg( i j - (i + j) \frac{\xi_1 \partial_1 Q(\xi) + \xi_2 \partial_2 Q(\xi)}{Q(\xi)} - \frac{\xi_1 \xi_2 \partial_1 \partial_2 Q(\xi)}{Q(\xi)} + 2 \frac{\xi_1 \partial_1 Q(\xi)}{Q(\xi)} \frac{\xi_2 \partial_2 Q(\xi)}{Q(\xi)} \bigg).
\end{equation}
On the other hand, we check that
$$\varepsilon^{\kappa_{i, j}} |\xi_1|^i |\xi_2|^j \leq 4 Q(\xi), \ |\xi_k| |\partial_k Q(\xi)| \leq 4 Q(\xi), \ {\rm and} \ |\xi_1| |\xi_2| |\partial_1 \partial_2 Q(\xi)| \leq 4 Q(\xi),$$
so that, by \eqref{namur} and \eqref{charleroi}, there exists some universal constant $K$ such that
$$\varepsilon^{\kappa_{i, j}} M \big( \widehat{K_\varepsilon^{i, j}} \big) \leq K.$$
Inequality \eqref{multiLq} then follows from \eqref{ibanez} applying Theorem \ref{thelizo}.
\end{proof}

\section{Sobolev bounds for $N_\p$ and $\Theta_\p$}
\label{Unity}

This section is devoted to the proof of the Sobolev estimates of $N_\p$, $\partial_1 \Theta_\p$ and $\partial_2 \Theta_\p$ stated in Proposition \ref{Sobolevbound}. As previously mentioned in Section \ref{Lent}, we focus on Sobolev bounds on $N_\p$.

\begin{prop}
\label{Recursive}
Let $\alpha \in \N^2$ and $1 < q < + \infty$. There exists some constant $K(q, \alpha)$, depending possibly on $\alpha$ and $q$, but not on $\p$, such that
\begin{equation}
\label{quirecure}
\begin{split}
& \| \partial^\alpha N_\p \|_{L^q(\R^2)} + \| \partial_1 \partial^\alpha N_\p \|_{L^q(\R^2)} + \| \partial_2 \partial^\alpha N_\p \|_{L^q(\R^2)}\\
+ & \| \partial_1^2 \partial^\alpha N_\p \|_{L^q(\R^2)} + \varepsilon_\p \| \partial_1 \partial_2 \partial^\alpha N_\p \|_{L^q(\R^2)} + \varepsilon_\p^2 \| \partial_2^2 \partial^\alpha N_\p \|_{L^q(\R^2)} \leq K(q, \alpha),
\end{split}
\end{equation}
for any $\p$ sufficiently small.
\end{prop}

\begin{remark}
The proof of Proposition \ref{Recursive} is by induction on the derivation order $\alpha$. The inductive assumption is given by \eqref{quirecure}. This explains the redundant form of this inequality.
\end{remark}

Proposition \ref{Sobolevbound} is a direct consequence of Proposition \ref{Recursive} invoking rescaled inequalities \eqref{thetap} and \eqref{thepat} to bound the functions $\partial_1 \Theta_\p$ and $\partial_2 \Theta_\p$.

\begin{proof}[Proof of Proposition \ref{Sobolevbound} (assuming Proposition \ref{Recursive})]
In view of \eqref{quirecure}, given any $k \in \N$ and any $1 < q < + \infty$, there exists some positive constant $K(k, q)$, not depending on $\p$, such that
\begin{equation}
\label{allez}
\| N_\p \|_{W^{k, q}(\R^2)} \leq K(k, q),
\end{equation}
for any $\p$ sufficiently small. In particular, by Sobolev embedding theorem,
\begin{equation}
\label{chardy}
\| N_\p \|_{\boC^k(\R^2)} \leq K(k).
\end{equation}
Using \eqref{allez} and \eqref{chardy}, inequality \eqref{thepat} becomes
\begin{equation}
\label{goodpat}
\Xi_\p(q, \alpha) \leq K(q, \alpha) \bigg( 1 + \varepsilon_\p^2 \sum_{0 \leq \beta < \alpha} \Xi_\p(q, \alpha - \beta) \bigg),
\end{equation}
where we have set as in Lemma \ref{Scaledtheta},
$$\Xi_\p(q, \alpha) \equiv \| \partial^\alpha \partial_1 \Theta_\p \|_{L^q(\R^2)} + \varepsilon_\p \| \partial^\alpha \partial_2 \Theta_\p \|_{L^q(\R^2)}.$$
By \eqref{thetap} and \eqref{allez}, the quantity $\Xi_\p(q, (0, 0))$ is bounded independently on $\p$, so that it follows by induction from formula \eqref{goodpat} that $\Xi_\p(q, \alpha)$ is bounded independently on $\p$ for any $1 < q < + \infty$ and any $\alpha \in \N^2$. Inequality \eqref{higher} follows invoking Sobolev embedding theorem for $q = + \infty$. This completes the proof of Proposition \ref{Sobolevbound}.
\end{proof}

The remainder of this section is devoted to the proof of Proposition \ref{Recursive}. As previously mentioned in Subsection \ref{Abricot}, the proof relies on decomposition \eqref{decomposition}. Recall that it is proved in \cite{Graveja3} that the functions $\eta$ and $\varphi$, and therefore $N_\p$ and $\Theta_\p$, belong to $W^{k, q}(\R^2)$ for any $k \in \N$ and any $1 < q \leq + \infty$. Hence, we can differentiate \eqref{decomposition} to any order $\alpha \in \N^2$ to obtain
\begin{equation}
\label{sapritch}
\partial^\alpha N_\p = K^{2, 0}_{\varepsilon_\p} \star \partial^\alpha f_\p + \varepsilon_\p^2 \sum_{i + j = 2} K_{\varepsilon_\p}^{i, j} \star \partial^\alpha \boR_{\varepsilon_\p}^{i, j}.
\end{equation}
Taking the $L^q$-norm of this expression and invoking the regularizing properties of the kernels provided by Lemma \ref{Multiprop}, we are led to
\begin{equation}
\label{alphaq0}
\| \partial^\alpha N_\p \|_{L^q(\R^2)} \leq K(q) \Big( \| \partial^\alpha f_\p \|_{L^q(\R^2)} + \varepsilon_\p^2 \sum_{i + j = 2} \| \partial^\alpha \boR_{\varepsilon_\p}^{i, j} \|_{L^q(\R^2)} \Big).
\end{equation}
In view of definitions \eqref{r20}, \eqref{r02}, \eqref{r11} and \eqref{effet}, the derivatives $\partial^\alpha f_\p$ and $\partial^\alpha \boR_{\varepsilon_\p}^{i, j}$ in the right-hand side of \eqref{alphaq0} are nonlinear functions of the derivatives of $N_\p$ and $\Theta_\p$, so that we may estimate their $L^q$-norms using Sobolev bounds on $N_\p$ and $\Theta_\p$.

This provides an iterative scheme to estimate the Sobolev norms of $N_\p$. Using the available information on the nonlinear source terms $f_\p$ and $\boR_{\varepsilon_\p}^{i, j}$, which is initially reduced to Lemma \ref{cosa}, we improve the regularity and integrability properties of $N_\p$ using inequality \eqref{alphaq0}. This in turn provides improved bounds of the nonlinear terms $f_\p$ and $\boR_{\varepsilon_\p}^{i, j}$.

As a consequence, we prove \eqref{quirecure} by induction on the derivation order $\alpha$. We first compute $L^q$-estimates of the nonlinear terms $f_\p$ and $\boR_{\varepsilon_\p}^{i, j}$, and of convolution equation \eqref{decomposition}. In particular, this requires to bound some derivatives of the phase $\Theta_\p$, which is made possible invoking Lemma \ref{Scaledtheta}. Using the initial bounds given by Lemma \ref{cosa}, we conclude that inequality \eqref{quirecure} holds for $\alpha = (0, 0)$. We then turn to higher order estimates. Assuming that \eqref{quirecure} holds for any index $\alpha$ such that $|\alpha| \leq k$, we derive $L^q$-estimates of the derivatives of order $k + 1$ of the functions $f_\p$ and $\boR_{\varepsilon_\p}^{i, j}$. In view of \eqref{alphaq0}, this provides bounds for the derivatives of order $k + 1$ of $N_\p$, so that we can prove that \eqref{quirecure} is also valid for any index $\alpha$ such that $|\alpha| = k + 1$. This completes the sketch of the proof of Proposition \ref{Recursive}, which is detailed below.

\subsection{$L^q$-estimates of nonlinear terms}
\label{nonlinear}

We first compute $L^q$-estimates on the nonlinear terms $f_\p$, $R_{\varepsilon_\p}^{i, j}$ and $\nu_{\varepsilon_\p}^{i, j}$.

\begin{lemma}
\label{Gaston}
Let $1 \leq q < + \infty$. There exists some universal constant $K$ such that
\begin{align}
\label{mosco}
& \| f_\p \|_{L^q(\R^2)} + \| R_{\varepsilon_\p}^{0, 2} \|_{L^q(\R^2)} + \varepsilon_\p \| R_{\varepsilon_\p}^{1, 1} \|_{L^q(\R^2)} \leq K \| N_\p \|_{L^{2 q}(\R^2)}^2,\\
\label{camba}
& \| R_{\varepsilon_\p}^{2, 0} \|_{L^q(\R^2)} \leq K \Big( \varepsilon_\p^{- 2} \| N_\p \|_{L^{2 q}(\R^2)}^2 + \| N_\p \|_{L^{3 q}(\R^2)}^3 + \| \partial_1 N_\p \|_{L^{2 q}(\R^2)}^2 \Big),\\
\label{sinclair}
& \| \nu_{\varepsilon_\p}^{2, 0} \|_{L^q(\R^2)} \leq K \Big( \varepsilon_\p^{- 2} \| N_\p \|_{L^{3 q}(\R^2)}^3 + \| \partial_2 N_\p \|_{L^{2 q}(\R^2)}^2 \Big),
\end{align}
and
\begin{equation}
\label{leguen}
\| \nu_{\varepsilon_\p}^{0, 2} \|_{L^q(\R^2)} \leq K \Big( \varepsilon_\p^{- 2} \| N_\p \|_{L^{2 q}(\R^2)}^2 + \| N_\p \|_{L^{3 q}(\R^2)}^3 + \| \partial_1 N_\p \|_{L^{2 q}(\R^2)}^2 + \varepsilon_\p^2 \| \partial_2 N_\p \|_{L^{2 q}(\R^2)}^2 \Big).
\end{equation}
\end{lemma}

\begin{proof}
Bounds \eqref{mosco}, \eqref{camba}, \eqref{sinclair} and \eqref{leguen} are consequences of inequalities \eqref{anquetil}, \eqref{poulidor} and \eqref{merckx} using H\"older inequalities. For the quantities involving the functions $\partial_1 \Theta_\p$ and $\partial_2 \Theta_\p$, we also use \eqref{thetap} to compute
$$\| (\partial_1 \Theta_\p)^2 \|_{L^q(\R^2)} + \varepsilon_\p \| N_\p \partial_2 \Theta_\p \|_{L^q(\R^2)} + \varepsilon_\p^2 \| (\partial_2 \Theta_\p)^2 \|_{L^q(\R^2)} \leq K(q) \| N_\p \|_{L^{2 q}(\R^2)}^2,$$
whereas
\begin{align*}
\| & N_\p (\partial_1 \Theta_\p)^2 \|_{L^q(\R^2)} + \varepsilon_\p^2 \| N_\p (\partial_2 \Theta_\p)^2 \|_{L^q(\R^2)}\\
\leq K(q) \| N_\p \|_{L^{3 q}(\R^2)} & \Big( \| \partial_1 \Theta_\p \|_{L^{3 q}(\R^2)}^2 + \varepsilon_\p^2 \| \partial_2 \Theta_\p \|_{L^{3 q}(\R^2)}^2 \Big) \leq K(q) \| N_\p \|_{L^{3 q}(\R^2)}^3.
\end{align*}
\end{proof}

\subsection{$L^q$-estimates of the convolution equation}
\label{Holder}

We now compute $L^q$-estimates of equation \eqref{decomposition} invoking the multiplier properties of the kernels $K_\varepsilon^{i, j}$ given by Lemma \ref{Multiprop}, and the previous $L^q$-estimates on the nonlinear terms $f_\p$, $R^{i, j}_{\varepsilon_\p}$ and $\nu_{\varepsilon_\p}^{i, j}$. This provides

\begin{lemma}
\label{thelq}
Let $1 < q < + \infty$. There exists some constant $K(q)$, depending only on $q$, such that
\begin{equation}
\label{keylq}
\begin{split}
& \| N_\p \|_{L^q(\R^2)} + \| \partial_1 N_\p \|_{L^q(\R^2)} + \| \partial_2 N_\p \|_{L^q(\R^2)}\\
+ \| \partial_1^2 & N_\p \|_{L^q(\R^2)} + \varepsilon_\p \| \partial_1 \partial_2 N_\p \|_{L^q(\R^2)} + \varepsilon_\p^2 \| \partial_2^2 N_\p \|_{L^q(\R^2)}\\
\leq K(q) \Big( \| N_\p & \|_{L^{2 q}(\R^2)}^2 + \varepsilon_\p^2 \| N_\p \|_{L^{3 q}(\R^2)}^3 + \varepsilon_\p^2 \| \partial_1 N_\p \|_{L^{2 q}(\R^2)}^2 + \varepsilon_\p^4 \| \partial_2 N_\p \|_{L^{2 q}(\R^2)}^2 \Big),
\end{split}
\end{equation}
for any $\p$ sufficiently small.
\end{lemma}

\begin{proof}
Given any $\alpha = (\alpha_1, \alpha_2)$ such that $0 \leq \alpha_1 + \alpha_2 \leq 2$, we estimate the $L^q$-norm of $\partial^\alpha N_\p$ using equations \eqref{decomposition}, so that
\begin{align*}
\| \partial^\alpha N_\p \|_{L^q(\R^2)} \leq & \| \partial^\alpha K_{\varepsilon_\p}^{2, 0} \star f_\p \|_{L^q(\R^2)} + \varepsilon_\p^2 \sum_{i + j = 2} \| \partial^\alpha K_{\varepsilon_\p}^{i, j} \star R_{\varepsilon_\p}^{i, j} \|_{L^q(\R^2)}\\
+ & \varepsilon_\p^4 \| \partial^\alpha K_{\varepsilon_\p}^{2, 0} \star \nu_{\varepsilon_\p}^{2, 0} \|_{L^q(\R^2)} + \varepsilon_\p^4 \| \partial^\alpha K_{\varepsilon_\p}^{0, 2} \star \nu_{\varepsilon_\p}^{0, 2} \|_{L^q(\R^2)}.
\end{align*}
Since by \eqref{defkernel},
$$\partial^\alpha K_{\varepsilon_\p}^{j, k} = i^{\alpha_1 + \alpha_2} K_{\varepsilon_\p}^{j + \alpha_1, k + \alpha_2},$$
the multiplier properties of Lemma \ref{Multiprop} provide
$$\| N_\p \|_{L^q(\R^2)} \leq K(q) \Big( \| f_\p \|_{L^q(\R^2)} + \varepsilon_\p^2 \sum_{i + j = 2} \| R_{\varepsilon_\p}^{i, j} \|_{L^q(\R^2)} + \varepsilon_\p^4 \| \nu_{\varepsilon_\p}^{2, 0} \|_{L^q(\R^2)} + \varepsilon_\p^4 \| \nu_{\varepsilon_\p}^{0, 2} \|_{L^q(\R^2)} \Big),$$
\begin{align*}
\| \partial_1 N_\p \|_{L^q(\R^2)} \leq & K(q) \Big( \| f_\p \|_{L^q(\R^2)} + \varepsilon_\p^2 \| R_{\varepsilon_\p}^{2, 0} \|_{L^q(\R^2)} + \varepsilon_\p^2 \| R_{\varepsilon_\p}^{1, 1} \|_{L^q(\R^2)} + \varepsilon_\p \| R_{\varepsilon_\p}^{0, 2} \|_{L^q(\R^2)}\\
& + \varepsilon_\p^4 \| \nu_{\varepsilon_\p}^{2, 0} \|_{L^q(\R^2)} + \varepsilon_\p^3 \| \nu_{\varepsilon_\p}^{0, 2} \|_{L^q(\R^2)} \Big),
\end{align*}
and
\begin{align*}
\| \partial_2 N_\p \|_{L^q(\R^2)} & + \| \partial_1^2 N_\p \|_{L^q(\R^2)} + \varepsilon_\p \| \partial_1 \partial_2 N_\p \|_{L^q(\R^2)} + \varepsilon_\p^2 \| \partial_2^2 N_\p \|_{L^q(\R^2)} \leq K(q) \bigg( \| f_\p \|_{L^q(\R^2)}\\
+ & \varepsilon_\p^2 \| R_{\varepsilon_\p}^{2, 0} \|_{L^q(\R^2)} + \varepsilon_\p \| R_{\varepsilon_\p}^{1, 1} \|_{L^q(\R^2)} + \| R_{\varepsilon_\p}^{0, 2} \|_{L^q(\R^2)} + \varepsilon_\p^4 \| \nu_{\varepsilon_\p}^{2, 0} \|_{L^q(\R^2)} + \varepsilon_\p^2 \| \nu_{\varepsilon_\p}^{0, 2} \|_{L^q(\R^2)} \bigg),
\end{align*}
Estimate \eqref{keylq} follows invoking nonlinear bounds \eqref{mosco}, \eqref{camba}, \eqref{sinclair} and \eqref{leguen}.
\end{proof}

\subsection{Initial bounds on $N_\p$ and its first order derivatives}
\label{initial}

In view of \eqref{keylq}, some preliminary $L^q$-bounds on $N_\p$, $\partial_1 N_\p$ and $\partial_2 N_ \p$ are required to inductively estimate the $L^q$-norms of these functions. These preliminary bounds are consequences of the uniform estimates given by \eqref{elinfiniestimates}, and the $L^2$-bounds provided by \eqref{primor}, \eqref{dial} and \eqref{badestimate1}.

\begin{lemma}
\label{Epidaure}
Let $2 \leq q \leq \frac{8}{3}$. There exists some constant $K(q)$, depending only on $q$, such that
\begin{equation}
\label{epi1}
\| N_\p \|_{L^q(\R^2)} \leq K(q),
\end{equation}
for any $\p$ sufficiently small. Moreover, given any $\frac{8}{3} < q < 8$, we have
\begin{equation}
\label{epi2}
\varepsilon_\p^\frac{2}{3} \| N_\p \|_{L^q(\R^2)} \leq K(q),
\end{equation}
whereas, given any $2 \leq q \leq + \infty$,
\begin{equation}
\label{epiderme}
\| \partial_1 N_\p \|_{L^q(\R^2)} + \varepsilon_\p \| \partial_2 N_\p \|_{L^q(\R^2)} \leq K(q) \varepsilon_\p^{\frac{6}{q} - 3}.
\end{equation}
\end{lemma}

\begin{proof}
For estimate \eqref{epiderme}, we have in view of \eqref{elinfiniestimates},
$$\| \partial_1 N_\p \|_{L^\infty(\R^2)} \leq \frac{K}{\varepsilon_\p^3}, \ {\rm and} \ \| \partial_2 N_\p \|_{L^\infty(\R^2)} \leq \frac{K}{\varepsilon_\p^4},$$
so that \eqref{epiderme} is a consequence of \eqref{dial} and \eqref{badestimate1} using standard interpolation between $L^q$-spaces.

The proofs of \eqref{epi1} and \eqref{epi2} are more involved. The first step is to compute $H^\alpha$-estimates of $N_\p$ combining equation \eqref{decomposition} with $H^\alpha$-bounds \eqref{crepewahou} on the kernels.

\begin{step}
\label{hachuncar}
Let $0 \leq \alpha \leq \frac{1}{4}$. There exists some constant $K(\alpha)$ such that
\begin{equation}
\label{hachuncar1}
\| N_\p \|_{H^\alpha(\R^2)} \leq K(\alpha),
\end{equation}
for any $\p$ sufficiently small. In particular, there exists some constant $K(q)$ such that \eqref{epi1} holds.
\end{step}

Applying Young inequality to decomposition \eqref{decomposition}, we have
\begin{align*}
& \| N_\p \|_ {H^\alpha(\R^2)} \leq \| K^{2, 0}_{\varepsilon_\p} \|_{H^\alpha(\R^2)} \Big( \| f_\p \|_{L^1(\R^2)} + \varepsilon_\p^2 \| R^{2, 0}_{\varepsilon_\p} \|_{L^1(\R^2)} + \varepsilon_\p^4 \| \nu^{2, 0}_{\varepsilon_\p} \|_{L^1(\R^2)} \Big)\\
+ & \varepsilon_\p^2 \| K^{1, 1}_{\varepsilon_\p} \|_{H^\alpha(\R^2)} \| \boR^{1, 1}_{\varepsilon_\p} \|_{L^1(\R^2)} + \varepsilon_\p^2 \| K^{0, 2}_{\varepsilon_\p} \|_{H^\alpha(\R^2)} \Big( \| R^{0, 2}_{\varepsilon_\p} \|_{L^1(\R^2)} + \varepsilon_\p^2 \| \nu^{0, 2}_{\varepsilon_\p} \|_{L^1(\R^2)} \Big).
\end{align*}
Combining \eqref{crepewahou} with \eqref{primor}, \eqref{cosa1} and \eqref{cosa2}, we derive \eqref{hachuncar1}, whereas \eqref{epi1} is a consequence of Sobolev embedding theorem,
$$H^\alpha(\R^2) \hookrightarrow L^q(\R^2),$$
for any $2 \leq q \leq \frac{2}{1 - \alpha}$.

The second step is to compute uniform bounds on $N_\p$ using Sobolev embedding theorem.

\begin{step}
\label{lymphatique}
Let $\nu > 0$. There exists some constant $K(\nu)$ such that
\begin{equation}
\label{linfini1}
\| N_\p \|_{L^\infty(\R^2)} \leq K(\nu) \Big( 1 + \varepsilon_\p^{- 1 - \nu} \Big),
\end{equation}
for any $\p$ sufficiently small.
\end{step}

In view of \eqref{epi1} and \eqref{epiderme}, there exists some number $q > 2$ such that
$$\| N_\p \|_{W^{1, q}(\R^2)} \leq K(\nu) \Big( 1 + \varepsilon_\p^{- 1 - \nu} \Big).$$
Estimate \eqref{linfini1} follows by Sobolev embedding theorem.

Combining with \eqref{epi1}, and invoking standard interpolation between $L^q$-spaces, estimate \eqref{linfini1} yields \eqref{epi2}.
\end{proof}

\subsection{Proof of inductive assumption \eqref{quirecure} for $\alpha = (0, 0)$}
\label{hypo0proved}

We now rely on Lemma \ref{thelq} to improve the preliminary estimates of Lemma \ref{Epidaure}. This gives

\begin{lemma}
\label{Union}
Let $1 < q < + \infty$. Then, assumption \eqref{quirecure} holds for $\alpha = (0, 0)$, i.e. there exists some constant $K(q)$, not depending on $\p$, such that
\begin{equation}
\label{grant}
\begin{split}
& \| N_\p \|_{L^q(\R^2)} + \| \partial_1 N_\p \|_{L^q(\R^2)} + \| \partial_2 N_\p \|_{L^q(\R^2)} \\
+ \| \partial_1^2 N_\p & \|_{L^q(\R^2)} + \varepsilon_\p \| \partial_1 \partial_2 N_\p \|_{L^q(\R^2)} + \varepsilon_\p^2 \| \partial_2^2 N_\p \|_{L^q(\R^2)} \leq K(q),
\end{split}
\end{equation}
for any $\p$ sufficiently small.
\end{lemma}

\begin{proof}
The proof relies on some bootstrap argument. Given any $1 < q \leq \frac{4}{3}$, we deduce from \eqref{keylq}, \eqref{epi1}, \eqref{epi2} and \eqref{epiderme},that
\begin{align*}
& \| N_\p \|_{L^q(\R^2)} + \| \partial_1 N_\p \|_{L^q(\R^2)} + \| \partial_2 N_\p \|_{L^q(\R^2)}\\
+ \| \partial_1^2 N_\p & \|_{L^q(\R^2)} + \varepsilon_\p \| \partial_1 \partial_2 N_\p \|_{L^q(\R^2)} + \varepsilon_\p^2 \| \partial_2^2 N_\p \|_{L^q(\R^2)} \leq K(q),
\end{align*}
so that by Sobolev embedding theorem,
\begin{align*}
\| N_\p \|_{L^q(\R^2)} + \varepsilon_\p \| \partial_1 N_\p \|_{L^q(\R^2)} + \varepsilon_\p^2 \| \partial_2 N_\p \|_{L^q(\R^2)} \leq K(q),
\end{align*}
for any $1 < q \leq 4$. Invoking \eqref{keylq} and \eqref{epi2} once more time, we are led to
\begin{align*}
& \| N_\p \|_{L^q(\R^2)} + \| \partial_1 N_\p \|_{L^q(\R^2)} + \| \partial_2 N_\p \|_{L^q(\R^2)}\\
+ \| \partial_1^2 N_\p & \|_{L^q(\R^2)} + \varepsilon_\p \| \partial_1 \partial_2 N_\p \|_{L^q(\R^2)} + \varepsilon_\p^2 \| \partial_2^2 N_\p \|_{L^q(\R^2)} \leq K(q),
\end{align*}
for any $1 < q \leq 2$. In particular, we have by Sobolev embedding theorem,
$$\| N_\p \|_{L^q(\R^2)} + \varepsilon_\p \| \partial_1 N_\p \|_{L^q(\R^2)} + \varepsilon_\p^2 \| \partial_2 N_\p \|_{L^q(\R^2)} \leq K(q),$$
for any $1 < q < + \infty$, so that \eqref{keylq} now yields \eqref{grant} for any $1 < q < + \infty$. This completes the proof of Lemma \ref{Union}.
\end{proof}

\subsection{Higher order estimates of the nonlinear terms $f_\p$ and $\boR_{\varepsilon_\p}^{i, j}$}
\label{alphaqnonlina}

We now assume that assumption \eqref{quirecure} holds for any $1 < q < + \infty$ and any $\alpha \in \N^2$ such that $|\alpha| \leq k$, and prove that it remains valid when $|\alpha| = k + 1$. Invoking again equation \eqref{decomposition}, we first derive improved Sobolev bounds on the nonlinear terms $f_\p$ and $\boR_{\varepsilon_\p}^{i, j}$. In view of definitions \eqref{r20}, \eqref{r02}, \eqref{r11} and \eqref{effet}, this requires to compute $L^q$-bounds on the derivatives of $\Theta_\p$. Hence, we show

\begin{lemma}
\label{Tetralogie}
Let $k \in \N$, and assume that \eqref{quirecure} holds for any $1 < q < + \infty$ and any $\alpha \in \N^2$ such that $|\alpha| \leq k$. Then, there exist some positive constants $K(q, \alpha)$, not depending on $\p$, such that
\begin{equation}
\label{walkyrie}
\| \partial^\alpha \partial_1 \Theta_\p \|_{L^q(\R^2)} + \varepsilon_\p \| \partial^\alpha \partial_2 \Theta_\p \|_{L^q(\R^2)} \leq K(q, \alpha),
\end{equation}
for any $1 < q < + \infty$, any $\alpha \in \N^2$ such that $|\alpha| \leq k + 1$, and any $\p$ sufficiently small.
\end{lemma}

\begin{proof}
Inequality \eqref{walkyrie} is a consequence of \eqref{thepat}. Applying Sobolev embedding theorem to assumption \eqref{quirecure}, we have
$$\| N_\p \|_{C^k(\R^2)} \leq K(k),$$
where $K(k)$ is some positive constant, not depending on $\p$. Therefore, given any $\alpha \in \N^2$ such that $|\alpha| \leq k + 1$, \eqref{thepat} may be written as
\begin{align*}
& \| \partial^\alpha \partial_1 \Theta_\p \|_{L^q(\R^2)} + \varepsilon_\p \| \partial^\alpha \partial_2 \Theta_\p \|_{L^q(\R^2)}\\
\leq K(q, \alpha) \bigg( \| \partial^\alpha N_\p \|_{L^q(\R^2)} & + \varepsilon_\p^2 \sum_{0 \leq \beta < \alpha} \Big( \| \partial^{\alpha - \beta} \partial_1 \Theta_\p \|_{L^q(\R^2)} + \varepsilon_\p \| \partial^{\alpha - \beta} \partial_2 \Theta_\p \|_{L^q(\R^2)} \Big) \bigg).
\end{align*}
Denoting
$$S_k^q = \sum_{|\alpha| \leq k + 1} \Big( \| \partial^\alpha \partial_1 \Theta_\p \|_{L^q(\R^2)} + \varepsilon_\p \| \partial^\alpha \partial_2 \Theta_\p \|_{L^q(\R^2)} \Big),$$
we deduce that
$$S_k^q \leq K(q, \alpha) \Big( \varepsilon_\p^2 S_k^q + \sum_{|\alpha| \leq k + 1} \| \partial^\alpha N_\p \|_{L^q(\R^2)} \Big).$$
Combined with assumption \eqref{quirecure}, this provides \eqref{walkyrie} for any $\p$ sufficiently small.
\end{proof}

We now turn to $L^q$-estimates of the functions $f_\p$ and $\boR_{\varepsilon_\p}^{i, j}$.

\begin{lemma}
\label{Aragorn}
Let $k \in \N$, and assume that \eqref{quirecure} holds for any $1 < q < + \infty$ and any $\alpha \in \N^2$ such that $|\alpha| \leq k$. Then, there exist some positive constants $K(q, \alpha)$, not depending on $\p$, such that
\begin{equation}
\label{frodon}
\| \partial^\alpha f_\p \|_{L^q(\R^2)} + \| \partial^\alpha \boR_{\varepsilon_\p}^{0, 2} \|_{L^q(\R^2)} + \varepsilon_\p \| \partial^\alpha \boR_{\varepsilon_\p}^{1, 1} \|_{L^q(\R^2)} + \varepsilon_\p^2 \| \partial^\alpha \boR_{\varepsilon_\p}^{2, 0} \|_{L^q(\R^2)}\leq K(q, \alpha),
\end{equation}
for any $1 < q < + \infty$, any $\alpha \in \N^2$ such that $|\alpha| \leq k + 1$, and any $\p$ sufficiently small.
\end{lemma}

\begin{proof}
Lemma \ref{Aragorn} is a consequence of assumption \eqref{quirecure}, and Lemma \ref{Tetralogie}. For instance, applying Leibniz formula to definition \eqref{effet}, we have
$$\big| \partial^\alpha f_\p \big| \leq K(\alpha) \sum_{0 \leq \beta \leq \alpha}
\Big( \big| \partial^\beta N_\p \big| \big| \partial^{\alpha - \beta} N_\p \big| + \big| \partial^\beta \partial_1 \Theta_\p \big| \big| \partial^{\alpha - \beta} \partial_1 \Theta_\p \big| \Big),$$
so that, by \eqref{quirecure}, \eqref{walkyrie}, and H\"older inequality,
$$\| \partial^\alpha f_\p \|_{L^q(\R^2)} \leq K(q, \alpha).$$
The proof is identical for the function $\boR_{\varepsilon_\p}^{1, 1}$, which verifies, in view of \eqref{r11} and Leibniz formula,
$$\big| \partial^\alpha \boR_{\varepsilon_\p}^{1, 1} \big| \leq K(\alpha) \sum_{0 \leq \beta \leq \alpha}
\big| \partial^\beta N_\p \big| \big| \partial^{\alpha - \beta} \partial_2 \Theta_\p \big|.$$
Similarly, for $\partial^\alpha \boR_{\varepsilon_\p}^{2, 0}$ and $\partial^\alpha \boR_{\varepsilon_\p}^{0, 2}$, it follows from \eqref{quirecure}, \eqref{walkyrie} and Leibniz formula, that
\begin{equation}
\label{boromir}
\begin{split}
& \| \partial^\alpha \boR_{\varepsilon_\p}^{0, 2} \|_{L^q(\R^2)} + \varepsilon_\p^2 \| \partial^\alpha \boR_{\varepsilon_\p}^{2, 0} \|_{L^q(\R^2)}\\
\leq K(q, \alpha) \bigg( 1 + \varepsilon_\p^2 & \Big\| \partial^\alpha \Big( \frac{(\partial_1 N_\p)^2}{1 - \frac{\varepsilon_\p^2}{6} N_\p} \Big) \Big\|_{L^q(\R^2)} + \varepsilon_\p^4 \Big\| \partial^\alpha \Big( \frac{(\partial_2 N_\p)^2}{1 - \frac{\varepsilon_\p^2}{6} N_\p} \Big) \Big\|_{L^q(\R^2)} \bigg),
\end{split}
\end{equation}
so that the proof of \eqref{frodon} reduces to estimate the $L^q$-norms in the left-hand side of \eqref{boromir}. In view of \eqref{quirecure}, we deduce from Sobolev embedding theorem that
\begin{equation}
\label{legolas}
\| \partial^\beta N_\p \|_{L^\infty(\R^2)} \leq K(\beta),
\end{equation}
for any $\beta \in \R^2$ such that $\beta \leq k$ and any $\p$ sufficiently small. When $|\alpha| \leq k$, the chain rule theorem combined with \eqref{quirecure} and \eqref{legolas} again provides estimates \eqref{frodon}. When $|\alpha| = k + 1$, this argument yields
$$\varepsilon_\p^2 \bigg\| \partial^\alpha \Big( \frac{(\partial_1 N_\p)^2}{1 - \frac{\varepsilon_\p^2}{6} N_\p} \Big) \bigg\|_{L^q(\R^2)} \leq K(q, \alpha) \Big( 1 + \varepsilon_\p^2 \big\| \partial^\alpha \partial_1 N_\p \big\|_{L^q(\R^2)} \Big) \leq K(q, \alpha),$$
and
$$\varepsilon_\p^4 \bigg\| \partial^\alpha \Big( \frac{(\partial_2 N_\p)^2}{1 - \frac{\varepsilon_\p^2}{6} N_\p} \Big) \bigg\|_{L^q(\R^2)} \leq K(q, \alpha) \Big( 1 + \varepsilon_\p^4 \big\| \partial^\alpha \partial_2 N_\p \big\|_{L^q(\R^2)} \Big) \leq K(q, \alpha),$$
where we have used the estimates in the second line of \eqref{quirecure} for the second inequalities. Combined with \eqref{boromir}, this completes the proof of inequality \eqref{frodon}.
\end{proof}

\subsection{Proof of Proposition \ref{Recursive}}
\label{Termine}

We are now in position to conclude the inductive proof of Proposition \ref{Recursive}.

\begin{proof}[Proof of Proposition \ref{Recursive}]
Given any $k \in \N$, we assume that \eqref{quirecure} holds for any $1 < q < + \infty$ and any $\alpha \in \N^2$ such that $|\alpha| \leq k$, and consider some index $\gamma \in \N^2$ such that $|\gamma| = k + 1$. Invoking equation \eqref{sapritch} and the kernel estimates of Lemma \ref{Multiprop}, we compute
\begin{equation}
\label{alphaq1}
\| \partial^\gamma \partial_1 N_\p \|_{L^q(\R^2)} \leq K(q) \Big( \| \partial^\gamma f_\p \|_{L^q(\R^2)} + \varepsilon_\p^2 \big( \| \partial^\gamma \boR_{\varepsilon_\p}^{2, 0} \|_{L^q(\R^2)} + \| \partial^\gamma \boR_{\varepsilon_\p}^{1, 1} \|_{L^q(\R^2)} \big) + \varepsilon_\p \| \partial^\gamma \boR_{\varepsilon_\p}^{0, 2} \|_{L^q(\R^2)} \Big),
\end{equation}
and
\begin{equation}
\label{alphaq2}
\begin{split}
& \| \partial^\gamma \partial_2 N_\p \|_{L^q(\R^2)} + \| \partial^\gamma \partial_1^2 N_\p \|_{L^q(\R^2)} + \varepsilon_\p \| \partial^\gamma \partial_1 \partial_2 N_\p \|_{L^q(\R^2)} + \varepsilon_\p^2 \| \partial^\gamma \partial_2^2 N_\p \|_{L^q(\R^2)}\\
\leq & K(q) \Big( \| \partial^\gamma f_\p \|_{L^q(\R^2)} + \varepsilon_\p^2 \| \partial^\gamma \boR_{\varepsilon_\p}^{2, 0} \|_{L^q(\R^2)} + \varepsilon_\p \| \partial^\gamma \boR_{\varepsilon_\p}^{1, 1} \|_{L^q(\R^2)} \big) + \| \partial^\gamma \boR_{\varepsilon_\p}^{0, 2} \|_{L^q(\R^2)} \Big).
\end{split}
\end{equation}
In view of inequalities \eqref{alphaq0}, \eqref{alphaq1} and \eqref{alphaq2}, and estimates \eqref{frodon}, assumption \eqref{quirecure} also holds for $\alpha = \gamma$. This completes the inductive proof of Proposition \ref{Recursive}.
\end{proof}

\section{Convergence towards \eqref{KP}}
\label{Concentration}

This section is devoted to the proofs of Theorem \ref{convGPKP} and Proposition \ref{bangkok}. As mentioned above in the introduction, our strategy is to prove that the sequence $(\partial_1 \Theta_\p)_{\p > 0}$ is, for $\p$ sufficiently small, a minimizing sequence for minimization problem \eqref{minikp} We then invoke Proposition \ref{ConcKP} to obtain the strong convergence of some subsequence towards a function $N_0$, which is a solution to minimization problem \eqref{minikp}, i.e. a ground state for \eqref{KP}. Finally, we improve the convergence using the previous Sobolev estimates.

\subsection{Weak convergence towards \eqref{KP}}
\label{faiblard}

We first use the Sobolev bounds provided by Proposition \ref{Sobolevbound} to establish the weak convergence of some subsequence $(N_{\p_n})_{n \in \N}$ to some non-constant solution $N_0$ to \eqref{SW}, as $\p_n \to 0$.

\begin{prop}
\label{Weekend}
There exists a subsequence $(\p_n)_{n \in \N}$, tending to $0$ as $n \to + \infty$, and a non-constant solution $N_0$ to \eqref{SW} such that, given any $1 < q < + \infty$,
\begin{equation}
\label{weaklq}
N_{\p_n} \rightharpoonup N_0 \ {\rm in} \ W^{1, q}(\R^2), \ {\rm as} \ n \to + \infty.
\end{equation}
In particular, given any $0 \leq \gamma < 1$, we have
\begin{equation}
\label{convholder}
N_{\p_n} \to N_0 \ {\rm in} \ \boC^{0, \gamma}(K), \ {\rm as} \ n \to + \infty,
\end{equation}
for any compact subset $K$ of $\R^2$.
\end{prop}

\begin{proof}
In view of bounds \eqref{higher}, there exists a subsequence $(\p_n)_{n \in \N}$, tending to $0$ as $n \to + \infty$, and a function $N_0$ such that \eqref{weaklq} holds for any $1 < q < + \infty$. Convergences \eqref{convholder} follow by standard compactness theorems. The proof of Proposition \ref{Weekend} therefore reduces to prove Lemma \ref{roti}, i.e. to establish that $N_0$ is a non-constant solution to \eqref{SW}.
\end{proof}

\begin{proof}[Proof of Lemma \ref{roti}]
Denoting
$$N_\p^0 = \frac{1}{2} K_{\varepsilon_\p}^{2, 0} \star f_\p,$$
we deduce from \eqref{decomposition} and Lemma \ref{noyaudur} that
$$\| N_\p - N_\p^0 \|_{L^2(\R^2)} \leq \varepsilon_\p^2 \sum_{i + j = 2} \| K_{\varepsilon_\p}^{i, j} \star \boR_{\varepsilon_\p}^{i, j} \|_{L^2(\R^2)} \leq \varepsilon_\p^2 \| \boR_{\varepsilon_\p}^{2, 0} \|_{L^1(\R^2)} + \varepsilon_\p^\frac{3}{2} \| \boR_{\varepsilon_\p}^{1, 1} \|_{L^1(\R^2)} + \varepsilon_\p^\frac{1}{2} \| \boR_{\varepsilon_\p}^{0, 2} \|_{L^1(\R^2)}.$$
In view of estimates \eqref{mosco}, \eqref{camba}, \eqref{sinclair} and \eqref{leguen}, and $L^q$-bounds \eqref{higher}, we obtain
$$\| N_\p - N_\p^0 \|_{L^2(\R^2)} \leq K \varepsilon_\p^\frac{1}{2},$$
so that
\begin{equation}
\label{restel2}
N_\p - N_\p^0 \to 0 \ {\rm in} \ L^2(\R^2), \ {\rm as} \ \p \to 0.
\end{equation}
We now claim that, up to some subsequence $(\p_n)_{n \in \N}$ satisfying \eqref{convholder},
\begin{equation}
\label{theclaim}
N_{\p_n}^0 \rightharpoonup \frac{1}{2} K_0 \star N_0^2 \ {\rm in} \ L^2(\R^2), \ {\rm as} \ n \to + \infty.
\end{equation}
Invoking the weak $L^2$-convergence provided by \eqref{weaklq}, we deduce from \eqref{restel2} and \eqref{theclaim} that the function $N_0$ satisfies
$$N_0 = \frac{1}{2} K_0 \star N_0^2,$$
so that, in view of \eqref{convSW}, the function $N_0$ is solution to \eqref{SW}.

Finally, in view of \eqref{pastriv} and convergences \eqref{convholder}, we have
$$N_0(0) \geq \frac{3}{5},$$
so that $N_0$ cannot be a constant solution to \eqref{SW}. This ends the proof of Lemma \ref{roti}.
\end{proof}

We now show Claim \eqref{theclaim}.

\begin{proof}[Proof of Claim \eqref{theclaim}]
Claim \eqref{theclaim} follows from \eqref{convholder} after the following simplification.

\setcounter{step}{0}
\begin{step}
\label{Reduc}
We have
$$N_\p^0 - \frac{1}{2} K_0 \star N_\p^2 \to 0 \ {\rm in} \ L^2(\R^2), \ {\rm as} \ \p \to 0.$$
\end{step}

In view of \eqref{effet}, we have
$$N_\p^0 - \frac{1}{2} K_0 \star N_\p^2 = \Big( K_{\varepsilon_\p}^{2, 0} - K_0 \Big) \star \Big( \frac{1}{3} N_\p^2 + \frac{1}{6} (\partial_1 \Theta_\p)^2 \Big) + \frac{1}{6} K_0 \star \Big( (\partial_1 \Theta_\p)^2 - N_\p^2 \Big),$$
so that, by Young inequality, and estimates \eqref{higher},
\begin{equation}
\label{restconv}
\Big\| N_\p^0 - \frac{1}{2} K_0 \star N_\p^2 \Big\|_{L^2(\R^2)} \leq K \Big( \| K_{\varepsilon_\p}^{2, 0} - K_0 \|_{L^2(\R^2)} + \| K_0 \|_{L^2(\R^2)} \| \partial_1 \Theta_\p - N_\p \|_{L^2(\R^2)} \Big).
\end{equation}
In view of definitions \eqref{kernel0} and \eqref{defkernel}, we have
$$\widehat{K_{\varepsilon_\p}^{2, 0}}(\xi) \to \widehat{K_0}(\xi), \ {\rm as} \ \p \to 0,$$
and
$$0 \leq \widehat{K_{\varepsilon_\p}^{2, 0}}(\xi) \leq \widehat{K_0}(\xi),$$
for any $\varepsilon_\p \geq 0$ and any $\xi \neq 0$. Since $K_0$ belongs to $L^2(\R^2)$ by Lemma \ref{noyaudur}, it follows from the dominated convergence theorem that
$$\int_{\R^2} \Big| \widehat{K_{\varepsilon_\p}^{2, 0}}(\xi) - \widehat{K_0}(\xi) \Big|^2 d\xi \to 0, \ {\rm as} \ \varepsilon_\p \to 0.$$
Hence, by Plancherel formula, the first term in the right-hand side of \eqref{restconv} tends to $0$, as $\p \to 0$, whereas the second term also tends to $0$ by \eqref{L2fort}. This completes the proof of Step \ref{Reduc}.

Invoking Step \ref{Reduc}, the proof of Claim \eqref{theclaim} reduces to

\begin{step}
\label{convweak}
Given some subsequence $(\p_n)_{n \in \N}$ such that \eqref{convholder} holds, we have
$$K_0 \star N_{\p_n}^2 \rightharpoonup K_0 \star N_0^2 \ {\rm in} \ L^2(\R^2), \ {\rm as} \ n \to + \infty.$$
\end{step}

First notice that, in view of \eqref{higher}, there exists some constant $K$, not depending on $n$, such that 
$$\| K_0 \star \big( N_{\p_n}^2 - N_0^2 \big) \|_{L^2(\R^2)} \leq \| K_0 \|_{L^2(\R^2)} \| N_{\p_n}^2 - N_0^2 \|_{L^1(\R^2)} \leq K,$$
so that by density of $\boC_c^\infty(\R^2)$ into $L^2(\R^2)$, the proof of Step \ref{convweak} reduces to prove that
\begin{equation}
\label{kerconv}
\int_{\R^2} \Big( K_0 \star \big( N_{\p_n}^2 - N_0^2 \big) \Big) \psi \to 0, \ {\rm as} \ n \to + \infty,
\end{equation}
for any function $\psi \in \boC_c^\infty(\R^2)$. Moreover, given any $\delta > 0$, the density of $\boC_c^\infty(\R^2)$ into $L^2(\R^2)$ also implies the existence of a function $\kappa_\delta \in \boC_c^\infty(\R^2)$ such that
$$\| K_0 - \kappa_\delta \|_{L^2(\R^2)} \leq \delta.$$
Given any function $\psi \in \boC_c^\infty(\R^2)$, this gives by Young inequality,
$$\bigg| \int_{\R^2} \Big( K_0 \star \big( N_{\p_n}^2 - N_0^2 \big) \Big) \psi \bigg| \leq \bigg| \int_{\R^2} \Big( \kappa_\delta \star \big( N_{\p_n}^2 - N_0^2 \big) \Big) \psi \bigg| + \delta \| N_{\p_n}^2 - N_0^2 \|_{L^1(\R^2)} \| \psi \|_{L^2(\R^2)},$$
which may be written as
$$\bigg| \int_{\R^2} \Big( K_0 \star \big( N_{\p_n}^2 - N_0^2 \big) \Big) \psi \bigg| \leq \bigg| \int_{\R^2} \big( \check{\kappa}_\delta \star \psi \big) \big( N_{\p_n}^2 - N_0^2 \big) \bigg| + K \delta,$$
denoting $\check{\kappa}_\delta(x) = \kappa_\delta(- x)$, and invoking \eqref{higher} and Fubini theorem. Since the function $\check{\kappa}_\delta \star \psi$ belongs to $\boC_c^\infty(\R^2)$, we deduce from \eqref{convholder} that
$$\int_{\R^2} \big( \check{\kappa}_\delta \star \psi \big) \big( N_{\p_n}^2 - N_0^2 \big) \to 0, \ {\rm as} \ n \to + \infty,$$
so that \eqref{kerconv} holds. This completes the proof of Step \ref{convweak} and of Claim \eqref{theclaim}.
\end{proof}

\subsection{Convergence of the energies}
\label{austin}

In order to apply Proposition \ref{Goodconc} to the family $(\partial_1 \Theta_\p)_{\p > 0}$ to deduce its strong convergence in the space $Y(\R^2)$, we first prove

\begin{prop}
\label{rover}
Let $(\p_n)_{n \geq 0}$ denote some subsequence, tending to $0$ as $n$ tends to $+ \infty$, such that \eqref{weaklq} and \eqref{convholder} hold. Then, up to some further subsequence, there exists a positive number $\mu_0$ such that
\begin{equation}
\label{convenergy}
E_{KP}(\partial_1 \Theta_{\p_n}) \to \boE_{\min}^{KP}(\mu_0), \ {\rm and} \ \int_{\R^2} |\partial_1 \Theta_{\p_n}|^2 \to \mu_0, \ {\rm as} \ n \to + \infty.
\end{equation}
\end{prop}

Proposition \ref{rover} is a consequence of Lemmas \ref{discrepancy} and \ref{upperbound}, so that we first address the proof of Lemma \ref{discrepancy}.

\begin{proof}[Proof of Lemma \ref{discrepancy}]
In view of formulae \eqref{p-slow} and \eqref{E-slow}, the discrepancy quantity $\Sigma(u_\p) = \sqrt{2} p(u_\p) - E(u_\p)$ may be recast in the slow space variables as
\begin{align*}
\Sigma(u_\p) = & - \sqrt{2} \frac{\varepsilon_\p}{144} \bigg( \int_{\R^2} \big( N_\p - \partial_1 \Theta_\p \big)^2 + \varepsilon_\p^2 \int_{\R^2} \Big( \frac{1}{2} (\partial_1 N_\p)^2 + \frac{1}{2} (\partial_2 \Theta_\p)^2 - \frac{1}{6} N_\p (\partial_1 \Theta_\p)^2 \Big)\\
+ & \varepsilon_\p^4 \int_{\R^2} \bigg( \frac{(\partial_2 N_\p)^2}{4 - \frac{4 \varepsilon_\p^2}{6} N_\p} + \frac{N_\p (\partial_1 N_\p)^2}{12 - 2 \varepsilon_\p^2 N_\p} - \frac{1}{12} N_\p (\partial_2 \Theta_\p)^2 \bigg) \bigg).
\end{align*}
Hence, we deduce from Proposition \ref{Sobolevbound} and estimate \eqref{dial} for the function $\partial_2 \Theta_\p$ that
\begin{equation}
\label{presquefini}
\Sigma(u_\p) = - \sqrt{2} \frac{\varepsilon_\p}{144} \bigg( \int_{\R^2} \big( N_\p - \partial_1 \Theta_\p \big)^2 + \varepsilon_\p^2 \int_{\R^2} \Big( \frac{1}{2} (\partial_1 N_\p)^2 + \frac{1}{2} (\partial_2 \Theta_\p)^2 - \frac{1}{6} N_\p (\partial_1 \Theta_\p)^2 \Big) + \underset{\p \to 0}{o} \big( \varepsilon_\p^2 \big) \bigg).
\end{equation}
Let us now recall that the value of $E_{KP}(\partial_1 \Theta_\p)$ is given by
$$E_{KP}(\partial_1 \Theta_\p) = \int_{\R^2} \Big( \frac{1}{2} (\partial_1^2 \Theta_\p)^2 + \frac{1}{2} (\partial_2 \Theta_\p)^2 - \frac{1}{6} (\partial_1 \Theta_\p)^3 \Big).$$
In particular, provided we may prove that
\begin{equation}
\label{ray1}
\| \partial_1 N_\p - \partial_1^2 \Theta_\p \|_{L^2(\R^2)} \to 0, \ {\rm as} \ \p \to 0,
\end{equation}
we have, in view of \eqref{higher} and \eqref{L2fort},
\begin{equation}
\label{sansnumero}
\int_{\R^2} \Big( \frac{1}{2} (\partial_1 N_\p)^2 + \frac{1}{2} (\partial_2 \Theta_\p)^2 - \frac{1}{6} N_\p (\partial_1 \Theta_\p)^2 \Big) - E_{KP}(\partial_1 \Theta_\p) \to 0, \ {\rm as} \ \p \to 0.
\end{equation}
Hence, by \eqref{presquefini},
\begin{equation}
\label{done}
\Sigma(u_\p) = - \sqrt{2} \frac{\varepsilon_\p}{144} \bigg( \int_{\R^2} \big( N_\p - \partial_1 \Theta_\p \big)^2 + \varepsilon_\p^2 E_{KP}(\partial_1 \Theta_\p) + \underset{\p \to 0}{o} \big( \varepsilon_\p^2 \big) \bigg).
\end{equation}
We then claim that
\begin{equation}
\label{ray2}
\frac{1}{\varepsilon_\p^2} \int_{\R^2} \big( N_\p - \partial_1 \Theta_\p \big)^2 \to 0, \ {\rm as} \ \p \to 0,
\end{equation}
which gives \eqref{sigma} using \eqref{done}.

In order to complete the proof of Lemma \ref{discrepancy}, it only remains to prove Claims \eqref{ray1} and \eqref{ray2}. For Claim \eqref{ray1}, we invoke equation \eqref{slow-2} and the Sobolev estimates of Proposition \ref{Sobolevbound}. Taking the $L^2$-norm of \eqref{slow-2}, we deduce from \eqref{higher} that
$$\| \partial_1 N_\p - \partial_1^2 \Theta_\p \|_{L^2(\R^2)} \leq K \varepsilon_\p,$$
where $K$ is some universal constant. Claim \eqref{ray1} follows taking the limit $\p \to 0$. Similarly, for Claim \ref{ray2}, we take the $L^2$-norm of equation \eqref{slow-1}, and obtain by \eqref{higher},
$$\| N_\p - \partial_1 \Theta_\p \|_{L^2(\R^2)} \leq K \varepsilon_\p^2,$$
so that
$$\frac{1}{\varepsilon_\p^2} \int_{\R^2} \big( N_\p - \partial_1 \Theta_\p \big)^2 \leq K \varepsilon_\p^2 \to 0, \ {\rm as} \ \p \to 0.$$
This concludes the proof of Lemma \ref{discrepancy}.
\end{proof}

\begin{remark}
\label{equivenergie}
Equivalence \eqref{equivalence2} is a consequence of inequality \eqref{sansnumero}, since it will be proved in the sequel that the quantity $E_{KP}(\partial_1 \Theta_\p)$ has a nonzero limit as $\p \to 0$.
\end{remark}

We now turn to the proof of Lemma \ref{upperbound}.

\begin{proof}[Proof of Lemma \ref{upperbound}]
Lemma \ref{upperbound} is a consequence of estimate \eqref{estimE2} of Theorem \ref{dim2}. Combining \eqref{estimE2} with \eqref{epsilonetp} and \eqref{sigma}, we obtain
$$E_{KP}(\partial_1 \Theta_\p) \leq - \frac{6912 \p^3}{\boS_{KP}^2 \varepsilon_\p^3} + \underset{\p \to 0}{o} \big( 1 \big),$$
so that by formula \eqref{p-slow},
$$E_{KP}(\partial_1 \Theta_\p) \leq - \frac{1}{54 \boS_{KP}^2} \bigg( \int_{\R^2} N_\p \partial_1 \Theta_\p \bigg)^3 + \underset{\p \to 0}{o} \big( 1 \big).$$
In view of \eqref{L2fort}, we have
$$E_{KP}(\partial_1 \Theta_\p) \leq - \frac{1}{54 \boS_{KP}^2} \bigg( \int_{\R^2} (\partial_1 \Theta_\p)^2 \bigg)^3 + \underset{\p \to 0}{o} \big( 1 \big).$$
On the other hand, it follows from Lemma \ref{MiniKP} that
$$E_{KP}(\partial_1 \Theta_\p) \geq \boE_{\min}^{KP} \bigg( \int_{\R^2} (\partial_1 \Theta_\p)^2\bigg) =- \frac{1}{54 \boS_{KP}^2} \bigg( \int_{\R^2} (\partial_1 \Theta_\p)^2 \bigg)^3,$$
which completes the proof of Lemma \ref{upperbound}.
\end{proof}

We finally deduce Proposition \ref{rover} from Lemma \ref{upperbound}.

\begin{proof}[Proof of Proposition \ref{rover}]
In view of \eqref{L2fort} and \eqref{weaklq}, we have
$$\underset{n \to + \infty}{\liminf} \int_{\R^2} \big( \partial_1 \Theta_{\p_n} \big)^2 \geq \int_{\R^2} N_0^2,$$
so that we may assume up to some further subsequence, that
\begin{equation}
\label{clapton}
\int_{\R^2} \big( \partial_1 \Theta_{\p_n} \big)^2 \to \mu_0, \ {\rm as} \ n \to + \infty,
\end{equation}
where
$$\mu_0 \geq \int_{\R^2} N_0^2 > 0.$$
Assertion \eqref{convenergy} is then a consequence of \eqref{page}, \eqref{clapton}, and formula \eqref{minKP} of $\boE_{\min}^{KP}$.
\end{proof}

\subsection{Strong convergence towards \eqref{KP}}
\label{costaud}

We now show Proposition \ref{Yconv}. i.e. the strong convergence of the family $(N_\p)_{\p > 0}$ in $L^2(\R^2)$ (up to some subsequence).

\begin{proof}[Proof of Proposition \ref{Yconv}]
In view of Proposition \ref{rover}, we may construct a subsequence $(\p_n)_{n \in \N}$, tending to $0$ as $n \to + \infty$, and some positive number $\mu_0$ such that
$$E_{KP}(\partial_1 \Theta_{\p_n}) \to \boE_{\min}^{KP}(\mu_0), \ {\rm and} \ \int_{\R^2} |\partial_1 \Theta_{\p_n}|^2 \to \mu_0, \ {\rm as} \ n \to + \infty.$$
By Proposition \ref{Goodconc}, up to some further subsequence, there exists some points $(a_n)_{n \in \N}$ and a ground state solution $N_0$ to \eqref{SWsigma}, with $\sigma = \frac{\mu_0^2}{(\mu^*)^2}$, such that
$$\partial_1 \Theta_{\p_n}(\cdot - a_n) \to N_0 \ {\rm in} \ Y(\R^2), \ {\rm as} \ n \to + \infty.$$
By \eqref{L2fort}, we are led to
\begin{equation}
\label{fortl2}
N_{\p_n}(\cdot - a_n) \to N_0 \ {\rm in} \ L^2(\R^2), \ {\rm as} \ n \to + \infty.
\end{equation}
Invoking Proposition \ref{Weekend} for the subsequence $(N_{\p_n}(\cdot - a_n))_{n \in \N}$, there exists a non-constant solution $\tilde{N}_0$ to \eqref{SW} such that weak convergences \eqref{weaklq} hold, up to some further subsequence. In particular, by \eqref{fortl2}, $N_0 = \tilde{N}_0$, so that $N_0$ is a ground state of speed $1$ of \eqref{KP}.

In order to complete the proof of Proposition \ref{Yconv}, it is now necessary to drop the invariance by translation, i.e. to prove that convergences in $Y(\R^2)$ and in $L^2(\R^2)$ also hold for the sequences $(\partial_1 \Theta_{\p_n})_{n \in \N}$, respectively $(N_{\p_n})_{n \in \N}$. Assuming first that, up to some further subsequence, there exists some number $a$ such that
$$a_n \to a, \ {\rm as} \ n \to + \infty,$$
we obtain that
$$\partial_1 \Theta_{\p_n} \to N_0(\cdot + a) \ {\rm in} \ Y(\R^2), \ {\rm and} \ N_{\p_n} \to N_0(\cdot + a) \ {\rm in} \ L^2(\R^2), \ {\rm as} \ n \to + \infty,$$
using the continuity of the map $a \mapsto \psi(\cdot - a)$ from $\R$ to any space $L^q(\R^2)$ (with $1 < q < + \infty$).
Since the function $x \mapsto N_0(x + a)$ is still a ground state of speed $1$ of \eqref{KP}, this completes the proof of Proposition \ref{Yconv}.

Hence, it remains to prove that the sequence $(a_n)_{n \in \N}$ contains some bounded subsequence. Assuming by contradiction that this is false, we may construct some subsequence, still denoted $(a_n)_{n \in \N}$, such that
\begin{equation}
\label{faux}
a_n \to + \infty, \ {\rm as} \ n \to + \infty.
\end{equation}
In view of \eqref{pastriv} and \eqref{higher}, there exists some positive number $\delta$, not depending on $n$, such that
$$\int_{B(0, 1)} N_{\p_n}^2 \geq 2 \delta,$$
for any $n$ sufficiently large. By \eqref{fortl2}, we also have
$$\int_{B(0, 1)} |N_0(x + a_n) - N_{\p_n}(x)|^2 dx \to 0, \ {\rm as} \ n \to + \infty,$$
so that
$$\int_{B(0, 1)} |N_0(x + a_n)|^2 dx \geq \delta,$$
for any $n$ sufficiently large. However, it is proved in \cite{Graveja7} that there exists some positive constant $K$ such that
$$N_0(x) \leq \frac{K}{1 + |x|^2}, \ \forall x \in \R^2,$$
so that
$$\frac{10 K}{1 + |a_n|^2} \geq \delta,$$
for any $n$ sufficiently large. This provides a contradiction to \eqref{faux} and completes the proof of Proposition \ref{Yconv}.
\end{proof}

\subsection{Proofs of Theorem \ref{convGPKP} and Proposition \ref{bangkok}}
\label{derdesders}

We finally conclude the proofs of our main theorems.

\begin{proof}[Proof of Theorem \ref{convGPKP}]
In view of Propositions \ref{Sobolevbound} and \ref{Yconv}, given any $k \in \N$ and any $1 < q < + \infty$, the family $(N_\p)_{\p > 0}$ is bounded, uniformly with respect to $\p$ small, in $W^{k, q}(\R^2)$, and converges, up to some subsequence, to some ground state $N_0$ of \eqref{KP} in the space $L^2(\R^2)$, as $\p \to 0$. Hence, by standard interpolation theorem, it actually converges to $N_0$ in $W^{k, q}(\R^2)$. This concludes the proof of Theorem \ref{convGPKP}.
\end{proof}

\begin{proof}[Proof of Proposition \ref{bangkok}]
The proof is identical to the proof of Theorem \ref{convGPKP}, considering the function $\partial_1 \Theta_\p$ instead of $N_\p$, and noticing that $Y(\R^2)$ continuously embeds into $L^2(\R^2)$.
\end{proof}

\begin{acknowledgement}
The first and second authors acknowledge support from the ANR project JC05-51279, "\'Equations de Gross-Pitaevskii, d'Euler, et ph\'enom\`enes de concentration", of the French Ministry of Research. The third author acknowledges support from the ANR project ANR-07-BLAN-0250, "\'Equations non lin\'eaires dispersives", of the French Ministry of Research.
\end{acknowledgement}

\bibliographystyle{plain}
\bibliography{Bibliogr}

\end{document}